\documentclass[11pt]{article}
\usepackage{amsmath,amssymb,amsfonts}
\usepackage[latin1]{inputenc}
\usepackage{epsfig}
\numberwithin{equation}{section}
\newtheorem{thm}{Theorem}[section]

\newtheorem{alem}[thm]{Lemma}
\newtheorem{aprop}[thm]{Proposition}
\newtheorem{acor}[thm]{Corollary}
\newtheorem{arem}[thm]{Remark}
\newenvironment{adem}[1][]%
   {\ \\ {\bf Proof #1. }}%
   {\hfill\mbox{\rule{2 true mm}{3 true mm}}}%\vskip 2 ex\noindent}
   {\ \\ {\bf Example #1. }}%
   {\hfill\mbox{\rule{2 true mm}{3 true mm}}}%\vskip 2 ex\noindent}
\newcommand{\R}{\mathbb{R}}
\newcommand{\1}{\mathbf{1}}
\newcommand{\Q}{\mathbb{Q}}
\renewcommand{\P}{\mathbb{P}}

\newcommand{\E}{{\mathbb E}}

\newcommand{\N}{{\mathbb N}}

 \title{A trajectorial interpretation of the dissipations of entropy and Fisher information for stochastic differential equations}

 \author{J.Fontbona\thanks{Universidad de Chile, DIM-CMM, UMI(2807) UCHILE-CNRS,  Casilla 170-3, Correo 3,
Santiago-Chile, e-mail:fontbona@dim.uchile.cl. Partially supported by
Fondecyt Grant 1110923,
BASAL-Conicyt and Millenium Nucleus NC120062.}\and B.Jourdain\thanks{Universit\'e  Paris-Est, CERMICS (ENPC), INRIA, 6-8 av Blaise Pascal, Cit\'e
     Descartes, Champs sur Marne, 77455 Marne-la-Vall\'ee Cedex 2, France -
     e-mail : jourdain@cermics.enpc.fr. Supported by ANR-09-BLAN-0216-01
MEGAS.}}

\begin{document}
 \maketitle
\begin{abstract}
 The
dissipation of general convex entropies for continuous time Markov
processes can be described in terms of backward martingales with respect to the tail filtration. The relative entropy is the expected value of a backward submartingale. In the case of
(non necessarily reversible) Markov diffusion processes, we  use
Girsanov theory to explicit the Doob-Meyer decomposition of this submartingale. We deduce a stochastic analogue of the well known entropy dissipation
formula,  which is  valid for general convex entropies, including the total
variation distance. Under additional regularity assumptions, and using It\^o's calculus and ideas of 
 Arnold, Carlen and Ju \cite{Arnoldcarlenju}, we obtain moreover a new Bakry Emery criterion which ensures exponential convergence of the entropy to $0$. This criterion is non-intrisic since it depends on the square root of the diffusion matrix, and cannot be written only in terms of the diffusion matrix itself. We provide examples where the classic Bakry Emery criterion fails, but our non-intrisic criterion applies without modifying the law of the diffusion process.


\smallskip

\end{abstract}

\section*{Introduction}

We are interested in the long-time behaviour of solutions to the
stochastic differential equation
\begin{equation}
   dX_t=\sigma(X_t)dW_t+b(X_t)dt\label{diffhom}
\end{equation}
where $b:\R^d\to\R^{d}$,
$\sigma:\R^d\to\R^{d\otimes d'}$ and $W=(W_t, t\geq 0)$
is a standard Brownian motion in $\R^{d'}$.

In case \eqref{diffhom} admits a reversible probability measure, the celebrated Bakry Emery curvature dimension criterion which involves the generator, the carr\'e du champs and the iterated carr\'e du champs  is a sufficient condition for this reversible measure  to satisfy a Poincar\'e inequality and a logarithmic Sobolev inequality. From these inequalities, one can respectively deduce exponential convergence to $0$ as $t\to\infty$ of the chi-square distance or the relative entropy between the marginal at time $t$ of the process and its reversible measure. These results have been extended to more general entropy functionals  (see for instance \cite{AMTU}).

In general, even when the stochastic differential equation \eqref{diffhom} admits an invariant probability measure, this measure might be not reversible. It is well known from both a probabilistic point of view \cite{HHS} and the point of view of partial differential equations \cite{Arnoldcarlenju} that  a contribution  in the drift term, antisymmetric with respect to the invariant measure, may accelerate convergence to this invariant measure as $t\to\infty$.

Throughout this paper,   we assume \begin{itemize}
\item[$H0)$]    $U:[0,\infty)\to\R$  is a convex function such that  $\inf
U>-\infty$,
\end{itemize} 
and we  consider the $U-$entropy of a probability measure
$p$ on a measurable space $(E,{\cal E })$, with respect to another probability measure $q$ on $(E,{\cal E })$, defined  by
\begin{equation*}
   H_U(p|q)=\begin{cases}
      \int_{\R^d} U \left(\frac{dp}{dq}(x)\right)dq(x)
 \ \mbox{ if
}p\ll q\\
+\infty\mbox{ otherwise}.
   \end{cases}
\end{equation*}
The particular cases $U(x)=\mathbf{1}_{x>0}x\ln(x)$ and $U(x)=(x-1)^2$ respectively  correspond to the usual entropy and the
$\chi^2$-distance. For $U(x)=|x-1|$, $H_U(p|q)$ coincides with the total variation distance  when $p\ll q$. Notice that $U$ is  continuous on $(0,+\infty)$ and  that $U(0)\geq \lim_{x\to 0^+}U(x)$.

The primal goal of this work is to recover, by arguments using It\^o's stochastic calculus, the results of \cite{Arnoldcarlenju} and  \cite{AMTU} about the long-time behaviour of  the $U$-entropy of the law of $X_t$ with respect to the invariant measure. 
Our approach is based on the following simple remark, valid for   an arbitrary (possibly non-homogeneous) continuous-time Markov process $(X_t:t\geq 0)$ with values in a measurable space $(E,{\cal E })$: 

\medskip

If
% $P_0$ and  $Q_0$ are two probability measures on $E$, and 
we denote \begin{itemize} 
\item[$\bullet$]   by $P_t$ and $Q_t  $ the time marginal laws of
$ X_t$  when the   initial laws are $P_0$ and $Q_0$, respectively, and 
\item[$\bullet$] by $(X^{P}_t)_{t\geq 0}$ and
$(X^{Q}_t)_{t\geq 0}$ realizations of the process $(X_t)$ with  $X_0$   respectively distributed according to $P_0$ and $Q_0,$
\end{itemize}
then,  as soon as $H_U(P_t|Q_t)<+\infty$ for some $t\geq 0$, one has   $P_s\ll Q_s$  for all $ s\geq t$ and the process
$$\left(U\left(\frac{dP_s}{dQ_s}(X^{Q}_s)\right)\right)_{s\geq t}$$ is  a backward ${\cal F}_s$-submartingale  with respect to the filtration ${\cal F}_s:=\sigma(X^Q_r,r\geq s)$.  In fact,  it is easily deduced from the Markov property that  if  $P_t\ll Q_t$ for some $t\geq 0$, then the law of $(X^{P}_r)_{r\geq t}$ is absolutely continuous with respect to the one of $(X^{Q}_r)_{r\geq t}$ and moreover,   $P_s\ll Q_s$ for all $ s\geq t$  with  $\left(\frac{dP_s}{dQ_s}(X^{Q}_s)\right)_{s\geq t}$ a backward martingale  with respect to the filtration ${\cal F}_s$.  Jensen's inequality ensures  that $t\mapsto H_U(P_t\vert Q_t)$ is non-increasing and implies the remark. 

The  convergence of the $U-$entropy
\begin{equation}
  H_U(P_s|Q_s)=\E\left(\left(U\left(\frac{dP_s}{dQ_s}(X^{Q}_s)\right)\right)\right)\underset{s\to \infty }{\longrightarrow} \E\left(U\left( \E\left(\frac{dP_t}{dQ_t}(X^{Q}_t)\bigg|\cap_{s\geq 0}{\cal F}_s\right)\right)\right)<\infty\label{limentrel}
\end{equation}
 is then deduced from the  a.s. convergence of $\frac{dP_s}{dQ_s}(X^{Q}_s)$ to $\E\left(\frac{dP_t}{dQ_t}(X^{Q}_t)\bigg|\cap_{s\geq 0}{\cal F}_s\right)$  (the fact that for $r\geq t$, $\frac{dP_r}{dQ_r}(X^{Q}_r)=0$ a.s. on the set $\left\{\E\left(\frac{dP_t}{dQ_t}(X^{Q}_t)\bigg|\cap_{s\geq 0}{\cal F}_s\right)=0\right\}$ permits to cope with the possible discontinuity of $U$ at $0$).

\medskip

The first section of the paper is dedicated to time-inhomogeneous Markov diffusions given by the stochastic differential equation\begin{equation}\label{diff}
dX_t=\sigma (t,X_t) d W_t+  b(t,X_t)dt
\end{equation}
where $b:\R_+\times\R^d\to\R^{d}$,
$\sigma:\R_+\times\R^d\to\R^{d\otimes d'}$. Under assumptions that guarantee that for both initial laws, the time-reversed processes are still diffusions, we  use
Girsanov theory to explicit the Doob-Meyer decomposition of the submartingale $\left(U(\frac{dP_s}{dQ_s}(X^{Q}_s))\right)_{s\geq t}$. In this way, we obtain
a stochastic analogue of the well known entropy dissipation
formula,  valid for general convex entropies (including total
variation). Taking expectations in this formula, we recover the well known fact  that the $U$-entropy dissipation is equal to the $U$-Fisher information.  The proofs of the main results of this section are given in Appendix  \ref{proofssecentdiss}.

  It should be noticed that the  idea of considering a trajectorial interpretation of entropy to obtain functional inequalities is not new, at least for reversible diffusions  (see e.g.  the work of Cattiaux \cite{C} whose results are nevertheless of quite different nature).  However,  even in the reversible case,  time reversal  of a diffusion starting out of equilibrium modifies the dynamics of the diffusion. The backward martingale approach   takes this fact into account and  moreover permits  the use of It\^o's calculus  under  less regularity than {\it a priori}  needed when working in the forward time direction.  Its interest  thus  goes beyond the treatment of non-reversible situations.

In the second section, we  further suppose that the stochastic differential equation is time-homogeneous (i.e. of the form \eqref{diffhom}) and that it admits an invariant probability distribution, that is chosen as the initial law $Q_0$. Under additional regularity assumptions, and using It\^o's calculus and some
ideas close to Arnold, Carlen and Ju \cite{Arnoldcarlenju}, we obtain a new Bakry Emery criterion which ensures exponential convergence of the $U$-Fischer information to $0$ and therefore exponential convergence of the $U$-entropy to $0$. In addition, under this criterion, the invariant measure satisfies a $U$-convex Sobolev inequality. This criterion is non-intrisic : it depends on the square root $\sigma$ of the diffusion matrix $a=\sigma\sigma^*$ and cannot be written only in terms of the diffusion matrix itself whereas, under mild regularity assumptions on $b$ and $a$, the law of $(X_t)_{t\geq 0}$ solving \eqref{diffhom} is characterized by the associated martingale problem only written in terms of $a$ and $b$. The main results of this section are proved in Appendix \ref{pdissipfish}. In Appendix \ref{compacj}, we point out that our approach allows us to recover the results and criterion provided in \cite{Arnoldcarlenju}. We also highlight the difference between the arguments leading to each of the two criteria. Additionally, we provide a combined criterion.

 Last, we provide in the third section  two examples where the classic Bakry Emery criterion fails, but our non-intrisic criterion ensures exponential convergence to equilibrium without modifying the law of the diffusion process.

As future work, we plan to investigate how to choose the square root $\sigma$ of the diffusion matrix in order to maximize the rate of exponential convergence to equilibrium given by our non-intrisic Bakry Emery criterion.

Throughout this work, we use the convention of summation over repeated indexes. 

{\bf Acknowledgements :} We thank Tony Leli\`evre (CERMICS) for pointing out to us the paper of Arnold,  Carlen and Ju \cite{Arnoldcarlenju} at an early stage of this research. We also thank Anton Arnold (TU Wien) for suggesting an improvement of our non-intrinsic Bakry-Emery criterion (see Remark \ref{impnibec} below). The first author last thanks the  hospitality and partial support of CERMICS.

\section{Entropy dissipation for diffusion processes}\label{secentdiss}
 {From} now on we assume that $(X_t,t\geq 0)$ is a Markov diffusion process,  solution to the stochastic differential equation
\begin{equation}\label{difft}
   dX_t=\sigma(t,X_t)dW_t+b(t,X_t)dt
\end{equation}
where $W=(W_t, t\geq 0)$
is a standard Brownian motion in $\R^{d'}$ and  $b:\R_+\times \R^d\to\R^{d}$,
$\sigma:\R_+\times\R^d\to\R^{d\otimes d'}$  are mesurable coefficients satisfying conditions that  will be specified below.

For $P_0$ and $Q_0$ two probability measures on $\R^d$, we now denote by $(X^P_t)_{t\geq 0}$ and $(X^Q_t)_{t\geq 0}$ two solutions of \eqref{difft} with $X_0$ respectively distributed according to $P_0$ and $Q_0$. For $t\geq 0$, the law of $X^P_t$ (resp. $X^Q_t$) is denoted by $P_t$ (resp. $Q_t$). 

    Our first goal is to  explicitly describe  the backward submartingale $U\left(\frac{dP_s}{dQ_s}(X^{Q}_s)\right)$ when $P_0\ll Q_0$ and, as a byproduct, the decrease of  its expectation $H_U(P_s|Q_s)$.  
    In a way, this backward-in-time  approach  to  entropy  is converse to F\"ollmer's  approach  to the study of   time reversal  of  diffusion processes  \cite{Foll}  (see  \cite{Follwa} for the  infinite dimensional case)  based on the   stability under time reversal  of the usual  pathwise entropy. The latter  corresponds to $U(r)=r\ln r$ in   Remark   \ref{Pathentro}  below.

 We fix 
 a finite time-horizon $T\in(0,+\infty)$  in order to work with standard (forward) filtrations  by time  reversal in $[0,T]$.  Let us introduce some notation:

\begin{itemize}
\item[$\bullet$]  
$\Q^T$ (resp. $\P^T$ ) will   denote the law of  the time reversed processes $(X^{Q}_{T-t})_{t\leq T}$ (resp.
$(X^{P}_{T-t})_{t\leq T}$) on  the canonical space $C([0,T],\R^d)$.
\item[$\bullet$] $(Y_t)_{t\leq T}$ stands from now on for 
the canonical process on $C([0,T],\R^d)$ and  ${\cal
G}_t=\sigma(Y_s,0\leq s\leq t)$ denotes its natural (complete, right continuous) filtration. 
\item[$\bullet$]  In all the sequel, $\E^T$ will denote the expectation under the law $\Q^T$.
\end{itemize}
Whenever  $P_0\ll Q_0$,  by the Markov property we have 
$\P^{T}\ll \Q^{T}$ with  $\frac{d\P^{T}}{d\Q^{T}}=\frac{dP_0}{dQ_0}(Y_T)$ and 
\begin{equation}\label{DtT}
D_t\stackrel{\rm def}{=}\frac{d\P^{T }}{d \Q^{T}}\bigg|_{{\cal G}_t}=
\frac{dP_{T-t}}{dQ_{T-t}}(Y_t), \quad 0\leq t\leq T
\end{equation}
is a $\Q^{T}-{\cal G}_t$ martingale. Moreover,  $ H_U(P_s|Q_s)<+\infty$ for $s\in [0,T]$  if and only if
$\left(U\left(D_t\right)\right)_{ 0\leq t\leq
T-s}$  is a uniformly integrable $\Q^{T}-{\cal G}_t$ submartingale, in which case one has
$$ H_U(P_t|Q_t)= \E^T \left(U\left(D_{T-t}\right)\right) \mbox{ for all }t\in [s,T].$$

\begin{arem}\label{Pathentro}
If  $ \mathbb{H}_U(\mathbb{P}_1|\mathbb{P}_2)$ denotes the pathwise
$U-$entropy of a  probability measure
$\mathbb{P}_1$ on $ C([0,T],\R^d)$ with respect to a second probability measure $\mathbb{P}_2$:
\begin{equation*}
   \mathbb{H}_U(\mathbb{P}_1|\mathbb{P}_2):=\begin{cases}
      \int_{C([0,T],\R^d)} U \left(\frac{d\mathbb{P}_1}{d\mathbb{P}_2}(w)\right)d\mathbb{P}_2(w)
&  \ \mbox{ if
}\mathbb{P}_1\ll \mathbb{P}_2, \\
+\infty & \mbox{ otherwise}, 
   \end{cases}
\end{equation*}
we easily  deduce that
$H_U(P_0|Q_0)=\mathbb{H}_U\left(law \left( X_t^P:\, 0\leq t \leq T\right) \bigg|  law \left(X_t^Q:\,  0\leq t \leq T\right)\right) =\mathbb{H}_U(\mathbb{P}^{T}|\mathbb{Q}^{T})$.
\end{arem}

In order to use It\^o calculus to obtain the explicit form of the Girsanov density   $D_t $ as a   $\Q^{T}-{\cal G}_t$  martingale, and then deduce the Doob-Meyer decomposition of the submartingale $U\left(D_t\right)$, we will  
assume that the Markov processes  $(X^{Q}_{T-t},t\leq T)$  and $(X^{P}_{T-t},t\leq T)$ are
 diffusion processes as well.  Conditions  ensuring this fact have
been studied  in F\"ollmer \cite{Foll},  in Hausmann and
Pardoux \cite{HP}, in Pardoux \cite{Pa} and in Millet {\it et. al}
\cite{MNS} among others, who in particular  provide the semimartingale
decomposition of $(X^{Q}_{T-t},t\leq T)$ in its filtration.
 We  recall   in Theorem \ref{TRMNS} below  the general results in \cite{MNS} in a slightly more restrictive setting. The following conditions are needed:

\begin{itemize}
\item[$H1)$] 
For each $T>0$, $\sup_{t\in[0,T]}(|b(t,0)|+|\sigma(t,0)|)<+\infty$ and  for every  $A>0$ there is a constant  $K_{T,A}>0$ such that
\begin{equation*}
   |b(t,x)-b(t,y)| +
\sum_{i=1}^{d'} |\sigma_{\bullet  i}(t,x)-\sigma_{\bullet  i}(t,y)|\leq K_{T,A} |x-y|, \,  \forall t\in [0,T],\;\forall x,y \in
B(0,A), 
\end{equation*}
where $\sigma_{\bullet  i}$ denotes the $i$-th column of the matrix $\sigma$ and $B(0,A)$ is the ball of radius $A>0$ centered at the origin in $\R^d$.  
Moreover,
\begin{itemize}
\item[$H1)'$] the constants $K_{T,A}$ do not depend on $A$,
 or
\item[$H1)''$]  for each $s\geq 0$, equation \eqref{difft} starting at time $s$ is strictly conservative,
and for any bounded open set $D\subset \R^d$,
$$\sup_{x\in D} \sup_{s\in [0,T]}
\E \left\{ \exp\left[\int_s^T \left[4 B_{s,t}(x)+8
\sum_j|A^j_{s,t}(x)|^2\right]dt\right]\right\}<\infty,$$ where
$B_{s,t}(x)=\left[\sum_{i,k=1}^d \partial_i
b_k((t,X_{s,t}(x))^2\right]^{\frac{1}{2}}, \
A_{s,t}^j(x)=\left[\sum_{i,k=1}^d \partial_i
\sigma_{kj}(t,X_{s,t}(x))^2\right]^{\frac{1}{2}}$ and $X_{s,t}(x)$
denotes the solution to \eqref{difft} starting from $x$ at time
$s<t$.
\end{itemize}

\item[$H2)_Q$]  For each $t>0$, the law $Q_t(dx)$ of $X_t^Q$ has a density $q_t(x)$ with
respect to Lebesgue measure.
\item[$H3)_Q$] Setting  $a_{ij}=(\sigma\sigma^*)_{ij}$, 
for each $i=1,\dots , d$  the distributional derivative $\partial_j (a_{ij}(t,x)q_t(x))$
 is a locally integrable function on $[0,T]\times\R^d$:
 $$\int_{0}^T\int_D | \partial_j (a_{ij}(t,x)q_t(x)) |dx
 dt<\infty \mbox{
for any bounded open set }D\subset \R^d.$$
\end{itemize}

For $(t,x)\in
[0,T]\times \R^d$ we write 
\begin{itemize}
\item $\bar{a}_{ij}(t,x):=a_{ij}(T-t,x),  i,j=1,\dots,d,$ 
\item $\bar{b}^{Q}_i(t,x):=-b_i(T-t,x)+\frac{\partial_j  (a_{ij}(T-t,x)q_{T-t}(x))
}{q_{T-t}(x)}$ \qquad  (with  the convention that the term involving $\frac{1}{q_{T-t}(x)}$ is zero when $q_{T-t}(x)$ is zero)
\end{itemize}
and notice that $\bar{b}^{Q}(t,x)$ is 
defined $dt\otimes dx$ a.e. on $[0,T]\times\R^d$ under assumption $H3)_Q$.

\begin{thm}\label{TRMNS} Assume that $H1)$ and $H2)_Q$  hold.
\begin{itemize}
\item[a)] Suppose moreover $H3)_Q$.  Then,  $\Q^{T}$ is a solution to the  martingale
problem:
 $$(MP)_Q : \quad M^f_t:=f(Y_t)-f(Y_0)-\int_0^t
\frac{1}{2} \bar{a}_{ij}(s,Y_s)\partial_{ij} f(Y_s)+
\bar{b}_i^Q(s,Y_s)\partial_i f(Y_s)ds,\; t\in [0,T]$$ is a continuous
martingale with respect to the filtration $({\cal G}_t)$ for all
$f\in C_0^{\infty}(\R^d)$.
\item[b)] Let $\tilde{b}:\R_+\times\R^d\to\R^{d}$ and
$\tilde{\sigma}:\R_+\times\R^d\to\R^{d\otimes d'}$ be measurable
functions  such that $\int^{T}_0\int_D |\tilde{a}_{ij}(t,x)|+
|\tilde{b}_i(t,x)  |q_{T-t}(x)dx dt<\infty$ for any bounded open set $D\subset \R^d$. Assume moreover that
$\Q^{T}$ is a solution to the martingale problem  with respect to
$({\cal G}_t)$ for the generator ${\cal L}_tf(x)= \frac{1}{2}
\tilde{a}_{ij}(t,x)\partial_{ij} f(x)+
\tilde{b}_i(t,x)\partial_i f (x)$. Then $H3)_Q$ holds, $\tilde{b}=\bar{b}$ and
$\tilde{a}=\bar{a}$.
\end{itemize}
\end{thm}

\begin{adem}
 According to Theorem 3.3 \cite{MNS}, under  $H1)$, $H2)_Q$ and $H3)_Q$,
 $(M^f_t)_{t\in[0,T)}$ is a continuous ${\cal G}_t$-martingale under
  ${\mathbb Q}^T$.  When $f$ is $C^\infty$ on $\R^d$ and vanishes outside  $B(0,A)$, we have  \begin{align}
   \E^T\bigg(\int_0^T&|\bar{b}_i^Q(s,Y_s)||\partial_i f(Y_s)|ds\bigg)\notag\\&\leq \sup_{B(0,A)}
   |\nabla f|\left(T\sup_{[0,T]\times B(0,A)}|b(s,x)|+
   \int_{[0,T]\times B(0,A)}\sum_{i=1}^{d}\left| \partial_j (a_{ij}(s,x)q_s(x))\right|dsdx\right)\label{defmart}
\end{align}
where the right-hand side is finite under $H1)$ and $H3)_Q$.
This implies that $\E^T(|M^f_T|)<+\infty$, and together with $H1)$, that $(M_t^f)_{t\in[0,T]}$ is a continuous ${\cal G}_t$-martingale under ${\mathbb Q}^T$. 
Part b) follows from Theorem 3.3 in \cite{MNS}.
\end{adem}

\medskip

Assume $H1)$, $H2)_P$, $H2)_Q$, $H3)_P$ and $H3)_Q$. Then, under  $(MP)_Q$ and $(MP)_P$, the process $Y_t$ is respectively a weak solution to the SDEs
\begin{equation}\label{diffrevQ}
dX_t= \bar{\sigma} (t,X_t) d \bar{W}_t + \bar{b}^{Q}(t,X_t)dt,\quad t\in
[0,T]
\end{equation}
and
\begin{equation*}
dX_t=\bar{\sigma} (t,X_t) d \tilde{W}_t+ \bar{b}^{P}(t,X_t)dt ,\quad t\in
[0,T] \,, 
\end{equation*}
 where $\bar{\sigma}(t,x)=\sigma(T-t,x)$ and $\bar{W}$ and $\tilde{W}$ are $d'$ dimensional Brownian motions in  possibly enlarged probability spaces. If for all $t>0$, $x\mapsto p_t(x)$ and $x\mapsto q_t(x)$ are strictly positive and differentiable, then
the difference between the  drift terms of the two equations  is given by
\begin{equation*}
\begin{split}
\bar{b}_i^P(t,x)-\bar{b}_i^Q(t,x)= & \bar{a}_{ij}(t,\cdot )\partial_j \ln p_{T-t}(x) -  \bar{a}_{ij}(t,\cdot )\partial_j \ln q_{T-t}(x)
 \\ 
= &  \bar{a}_{i j}(t,x)\partial_j \left[ \ln\frac{  p_{T-t}}{  q_{T-t}}(x) \right]. \\
\end{split}
\end{equation*} 
If uniqueness in law holds for the second stochastic differential equation, then the simplest form of Girsanov theorem allows us to deduce that 
$$D_t=\frac{p_T}{q_T}(Y_0)\exp\bigg\{\int_0^t \nabla^* \left[ \ln\frac{  p_{T-t}}{  q_{T-t}}(Y_t) \right]\bar{\sigma} (t,Y_t) d \bar{W}_t -\frac{1}{2}\int_0^t \nabla^* \left[\ln
\frac{p_{T-s}}{q_{T-s}}\right](Y_s)  \bar{a}(s,Y_s) \nabla
\left[\ln
\frac{p_{T-s}}{q_{T-s}}\right](Y_s)ds\bigg\} $$
(in the above equation and from now on, we denote by $\nabla^*$ the transpose of the gradient ).
However,  in the general case when $q_t(x)$ or $p_t(x)$ may vanish  and are possibly not differentiable, it is not clear what sense should be given to the derivatives above. If the diffusion matrix is singular, neither is it clear that the difference of drift terms $\bar{b}^Q$ and $\bar{b}^P$ (defined by means of distributional derivatives) is in the range of the diffusion matrix, which is required in order  to  use Girsanov theorem. 

The problem of finding $D_t$ in the general case is reminiscent and, somehow, reciprocal to the  stochastic construction of Nelson processes, where $\Q^T$ and the possibly singular difference of drift terms are given, and one aims to construct $\P^T$ (see for instance \cite{CL}). The following technical lemma answers the question in the  most general situations covered by Theorem \ref{TRMNS}. Its proof, not hard but lengthy, relies on Girsanov theory in the absolutely continuous setting and is  given in the Appendix  \ref{pDtTexplicit} section.  Recall that an
element $\mathbb{P}_0\in {\cal M}$ of  a given  set ${\cal M}$ of
probability measures in $C([0,T],\R^d)$ is said to be {\it
extremal} if $\mathbb{P}_0=\alpha \mathbb{P}_1 + (1-\alpha)
\mathbb{P}_2$ for some $\mathbb{P}_1,\mathbb{P}_2\in{\cal M}$ and
$\alpha \in (0,1)$ implies $\mathbb{P}_0=\mathbb{P}_1 =
\mathbb{P}_2$.

\begin{alem}\label{DtTexplicit}
Assume  that $H1)$ , $H2)_Q$,  $H3)_Q$  and $H3)_P$ hold, with $P_0\ll Q_0$, and let  $\frac{p_t}{q_t}(x)$ be the Radon-Nikodyn derivative of   $ p_{t}(x)dx dt $ w.r.t   $ q_{t}(x)dx dt $  on $[0,T]\times \R^d$. 
 Then, 
\begin{itemize}
\item[a)] there exists a measurable function  in $[0,T]\times \R^d \to \R^d$ denoted $(t,x)\mapsto 
\nabla\ln [\frac{p_{t}}{q_{t}}](x)$ 
%and $(t,x)\mapsto  \nabla[\frac{p_{t}}{q_{t}}](x)$
 such that 
\begin{equation*}
\bar{b}^P(t,x)-\bar{b}^Q(t,x)=   \bar{a}(t,x)\nabla \left[ \ln\frac{  p_{T-t}}{  q_{T-t}}(x) \right] , \qquad p_{T-t}(x)dx \ dt \, a.e. .
\end{equation*} 

%\begin{equation*}
%\Q^T a.s.: \, 
%\int_0^T \left(\nabla
%\left[\frac{p_{T-s}}{q_{T-s}}\right](Y_s)\right)^* \bar{a}(s,Y_s)
%\nabla \left[\frac{p_{T-s}}{q_{T-s}}\right](Y_s) ds<\infty.
%\end{equation*}

\item[b)] Define $ q_{t}(x)dx \ dt$  a.e. in $[0,T]\times \R^d$ the function  $(t,x)\mapsto 
\nabla[\frac{p_{t}}{q_{t}}](x)$ by 
 \begin{equation*}
  \nabla \left[ \frac{p_{t}}{q_{t}}\right](x) : = \frac{p_{t}}{q_{t}}(x)   \nabla  \left[  \ln \frac{p_{t}}{q_{t}}\right](x) 
  \end{equation*} 
and assume  moreover that $\Q^T$ is an extremal solution to the
martingale problem $(MP)_Q$. Then, the $\Q^T$--$({\cal G}_t)$   martingale     $(D_t)_{t\in
[0,T]}$ introduced in \eqref{DtT}  has a continuous version (denoted in the same way) satisfying 
\begin{equation*}
\begin{split}
D_t= & \frac{ p_T}{ q_T}(Y_0) + \int_0^{t}D_s \nabla \ln
\left[\frac{p_{T-s}}{q_{T-s}}\right](Y_s)  \mathbf{1}_{s<R}  \cdot dM_ s\ \\
= &  \frac{ p_T}{q_T}(Y_0) + \int_0^{t} \nabla
\left[\frac{p_{T-s}}{q_{T-s}}\right](Y_s) \mathbf{1}_{\{
\frac{p_{T-s}}{q_{T-s}}(Y_s)>0\}}  \cdot
dM_ s \\
\end{split}
\end{equation*}
where  $M_t=(M_t^i)_{i=1}^d$ are the  continuous local
martingales w.r.t.  $\Q^T$ and $({\cal G}_t)$ defined by 
$$M_t^i:=Y_t^i-Y_0^i-\int_0^{t} \bar{b}_i^Q(s,Y_s)ds, \ t\in [0,T]$$
 and
 $R$ is  the $({\cal G}_t)$-stopping time $R:=\inf\{s
\in [0,T]:D_s=0\}$.
Moreover, $\Q^T$ a.s.,   one has
$$\langle D \rangle_t=\int_0^t \nabla^*
\left[\frac{p_{T-s}}{q_{T-s}}\right](Y_s)  \bar{a}(s,Y_s)
\nabla \left[\frac{p_{T-s}}{q_{T-s}}\right](Y_s)\mathbf{1}_{s<R} \ ds \,, \,  \forall t\in [0,T].$$
\end{itemize}
\end{alem}

{F}rom the proof of Lemma \ref{DtTexplicit}  it will be clear that if   $p_t$ and $q_t$ are everywhere strictly positive and of class $C^1$,   $(t,x)\mapsto 
\nabla [\frac{p_{t}}{q_{t}}](x)$ and $(t,x)\mapsto 
\nabla\ln [\frac{p_{t}}{q_{t}}](x)$   can be respectively  taken to be the usual gradient and  gradient  of the   logarithm of $\frac{p_{t}}{q_{t}}$. % An exponential form for $D_t$ can also be given. 

We now introduce the notations $U'_-$ and $U''(dy)$ for the left-hand derivative
of the restriction of the convex function $U:[0,\infty)\to \R$  to $(0,+\infty)$ and the non-negative measure on $(0,+\infty)$
equal to the second order distribution derivative of this restriction. 

We are ready to state  the main result of this section:

\begin{thm}\label{entprod}   (Stochastic $U$-Entropy dissipation)  Let $Q_0$ and $P_0$ be probability measures on $\R^d$ such that
$$H_U(P_0|Q_0)<\infty$$
and assume  that $H1)$ , $H2)_Q$,  $H3)_Q$  and $H3)_P$ hold. 
Suppose
moreover that $\Q^T$ is an extremal solution to the
martingale problem $(MP)_Q$.

 Then, the submartingale $(U(D_t))_{t\in [0,T]}$ has the Doob-Meyer decomposition
\begin{equation}\label{itoU}
\begin{split}
\forall t\in[0,T],\;U(D_t) =& U(D_0)+  \int_0^{t} U'_-(D_s) \nabla
\left[\frac{p_{T-s}}{q_{T-s}}\right](Y_s)\mathbf{1}_{s<R} \cdot
dM_ s
\\& +\frac{1}{2}\int_{(0,+\infty)} L^r_{t}(D)U''(dr) -1_{\{0<R\leq t\}}\Delta U(0),
\end{split}
\end{equation}
where  $R:=\inf\{s \in [0,T]:D_s=0\}$, $\Delta U(0)=\lim_{x\to
0^+}U(x)-U(0)\leq 0$ and $L^r_t(D)$ denotes the local time at
level $r\geq 0$ and time $t$ of  the continuous version of the martingale
$(D_s)_{s\in[0,T]}$.

In particular, if $U$ is continuous on $[0,+\infty)$ and $C^2$ on $(0,+\infty)$, one has
\begin{equation}\label{itoUC2}
\begin{split}
\forall t\in[0,T],\;U(D_t) =& U(D_0)+  \int_0^{t} U'(D_s) \nabla
\left[\frac{p_{T-s}}{q_{T-s}}\right](Y_s)\mathbf{1}_{s<R} \cdot
dM_ s
\\& +\frac{1}{2} \int_0^t
U''\left(\frac{p_{T-s}}{q_{T-s}}(Y_s)\right)\left(\nabla^*\left[
\frac{p_{T-s}}{q_{T-s}}\right] \bar{a}(s,\cdot) \nabla \left[
\frac{p_{T-s}}{q_{T-s}}\right]\right)(Y_s) \mathbf{1}_{s<R} ds. 
\end{split}
\end{equation}
\end{thm}

Theorem \ref{entprod}   is proved  in the Appendix \ref{proofssecentdiss} section. 
   We next briefly discuss some of   its assumptions and then state  some consequences.

\begin{arem}\label{condthmprincip}
\begin{itemize}
\item[a)] 
 By Theorem  3.1 in \cite{HP},   conditions $H2)_Q$ and $H3)_Q$  hold  under condition  $H1)'$ if
$Q_0$ has a density $q_0$ w.r.t. the Lebesgue measure
  s.t. $\int_{\R^d}\frac{q_0^2(x)dx}{1+|x|^k}<+\infty$ for some $k>0$ and either
 \begin{equation*}
\forall T>0,\;\exists\varepsilon>0,\;\forall
  (t,x)\in[0,T]\times\R^d,\;a(t,x){=}\sigma\sigma^*(t,x)\geq\varepsilon
  I_d,
\end{equation*}   or   the second order
  distribution derivatives
  $\frac{\partial^2a^{ij}}{\partial x_i\partial x_j}(t,x)$ are bounded
  on $[0,T]\times\R^d$ for each $T>0$ (by Theorem 3.1. in p. 1199   \cite{HP}, the latter conditions imply that
    (A)(ii) in  p. 1189 and thus Theorem 2.1 therein  hold). 
  In particular, under $H1)'$ and  the previous conditions, $H2)_P$ and $H3)_P$ also hold if  for instance $P_0\ll Q_0$ and $\frac{d P_0}{d Q_0}$ has polynomial growth.
  \item[b)]    Condition H1)''  introduced in  \cite{MNS}    allows us to include in our study the   fundamental examples of Langevin diffusions with  $a(x)=I_d$ and   $b(x)=-\nabla V(x)$ for  a nonnegative $C^2$ potential $V$, possibly superquadratic  but satisfying:
\begin{equation}
   \limsup_{|x|\to\infty}\frac{-x^*\nabla V(x)}{|x|^2}<+\infty,\;\limsup_{|x|\to\infty}\frac{\Delta V}{|\nabla V|^2}(x)<2\mbox{ and }\limsup_{|x|\to\infty}\frac{\sqrt{\partial_{ik}V\partial_{ik}V}}{V}(x)=0.\label{condv}
\end{equation}
See the Appendix section \ref{proofsuperquad} for a proof of this fact. 
\item[c)] Extremality of the solution $\Q^T$      to  the martingale problem   $(MP)_Q$  is implied by 
  pathwise uniqueness for the stochastic differential
equation \eqref{diffrevQ}. In the relevant case when $\sigma$ and $b$ in   \eqref{difft}  are  time-homogeneous and \eqref{diffhom} admits  an
invariant density   $p_{\infty}(x)>0$,  for the  choice  $Q_0(dx)=p_{\infty}(x)dx$  equation \eqref{diffrevQ}      takes the form 
\begin{equation*}
dX_t= \sigma (X_t) d W_t+\left(\frac{\partial_j(a_{\bullet j}p_\infty)}{p_\infty}(X_t)-b(X_t)\right)dt \quad t\in
[0,T].
\end{equation*}
 Pathwise uniqueness  for this SDE 
can be proved under H1) by a  standard argument using localization,  It\^o's formula and Gronwall's lemma,
 whenever the function
$- \frac{\partial_j
\left(a_{\bullet j}p_{\infty}\right)}{p_{\infty}}$ is the sum of a locally Lipschitz function and a monotone
function.  This is for instance the case when  $a =I_d$  and
$p_{\infty}(x)=C e^{-2V(x)}$ for some convex function $V:\R^d\to
\R$, or when  the strictly positive density $p_{\infty}$  and $a$  have locally Lipschitz derivatives.
 
  \end{itemize}
    \end{arem}

The proof of Theorem \ref{entprod} will justify that expectations can be taken in  \eqref{itoU} and \eqref{itoUC2},  yielding
\begin{acor}\label{expectentprod}  ($U$-Entropy dissipation)  
Under the assumptions of Theorem \ref{entprod},   
\begin{equation}
\forall t\in[0,T],\quad    H_U(P_t|Q_t)=H_U(P_T|Q_T)-\Delta U(0){\mathbb Q}^T(0<R\leq T-t) +\frac{1}{2} \int_{(0,+\infty)}\E^T\left( L^r_{T-t}(D) \right)U''(dr).  \label{entderivloc}
\end{equation}
% By taking $t=0$ and  then replacing in the obtained identity $T$ by $t$ we deduce
% \begin{equation}
%  \forall t\in [0,T], \quad     H_U(P_t|Q_t)=H_U(P_0|Q_0)+\Delta U(0){\mathbb Q}^t(0<R^{(t)}\leq t )-\frac{1}{2} \int_{(0,+\infty)}\E^t\left( L^r_{t}(D^{(t)}) \right) U''(dr),\label{entderivloc'}
% \end{equation}
% where for each  $t\in [0,T] $, $D^{(t)}$ stands for the process defined as in $\eqref{DtT}$ but with $t$  taking the role of $T$ therein, and $R^{(t)}$  is defined  as $R$ accordingly.

If  $U$ is moreover continuous on $[0,+\infty)$ and $C^2$ on $(0,+\infty)$,  we get the well known expression for the entropy dissipation: \begin{multline}\label{entderiv}
 \forall t\in[0,T] , \quad  H_U(P_t|Q_t)=H_U(P_0|Q_0) \\
-\frac{1}{2}\int_0^t
\int_{\{\frac{p_{s}}{q_{s}}(x)>0\}}U''\left(\frac{p_{s}}{q_{s}}(x)\right)\left(\nabla^*\left[
\frac{p_{s}}{q_{s}}\right] a(s,\cdot) \nabla \left[
\frac{p_{s}}{q_{s}}\right]\right)(x) q_{s}(x)dx ds,
\end{multline}
 with $U''(r)$  now  standing for the second order derivative of $U$ at  $r>0$.

  \end{acor}

The  particular  case  $U(x)=|x-1|$ of the  total variation  distance is more intricate but we are still able  to derive an analogous dissipation formula. To our knowledge this formula  is new: 

\begin{acor}\label{TV} (Dissipation of total variation)  Under the assumptions of Theorem \ref{entprod}, suppose  moreover  that  for a.e. $t\in [0, T]$, the functions $x\mapsto q_t(x)$ and $x\mapsto \frac{p_t}{q_t}(x)$ are respectively of class $C^1$ and $C^{2}$ and there exists a sequence $(r_n)_n$ of positive numbers tending to $+\infty$ as $n\to\infty$, such that
$\lim_{n\to \infty} \frac{1}{r_n}\int_{\{r_n\leq |x|<  2r_n \}} 
\left| a(t,x)
\nabla \left[\frac{p_t}{q_t}\right] (x)\right| q_t (x) dx  =0 $. Furthermore, assume that
$\int_0^{T}   \int_{\R^d}  \left| \nabla\cdot \left[ \bar{a}(s,x)
\nabla \left[\frac{p_{T-s}}{q_{T-s}}\right](x) q_{T-s}(x)\right] \right | dx ds <\infty $. Then,   $\forall t\in  [0,T]$, 
 $$\|P_t-Q_t\|_{\rm TV}=\|P_0-Q_0\|_{\rm TV} + \frac{1}{2}\int_0^{t}   \int_{\R^d}  \widetilde{sign} \left(\frac{p_{s}}{q_{s}}-1\right) (x)  \nabla\cdot \left[ a(s,x)
\nabla \left[\frac{p_{s}}{q_{s}}\right](x) q_{s}(x)\right] dx ds $$
 where $\widetilde{sign}(r)=- \mathbf{1}_{(-\infty,0)}(r)+  \mathbf{1}_{(0, \infty)}(r)$ and the integral is non-positive for all $t\in  [0,T]$.

\end{acor}

The proof is given in Appendix \ref{apptv}.

\begin{arem}
\begin{itemize}
\item[a)]  Denote by  $\Q$  the law of  $(X^Q_t,t\leq T)$ and by $\E$ the corresponding expectation. The following ``forward'' version of formula \eqref{entderivloc} 
holds under the assumptions of Theorem \ref{entprod} if moreover $\frac{p_t}{q_t}(Y_t)$ is a continuous $({\cal G}_t)$ semimartingale under $\Q$ (in particular if $(t,x)\mapsto \frac{dP_t}{dQ_t}(x)$ has a version of  class $C^{1,2}$):
\begin{equation*}
   \forall t\in  [0,T],  \quad  H_U(P_t|Q_t)=H_U(P_0|Q_0)+\Delta U(0){\mathbb Q}(0<S\leq t)\notag- \frac{1}{2}\int_{(0,+\infty)}\E\left( L^r_t \left(\frac{p_.}{q_.}(Y_.)\right)   \right) U''(dr),
   \end{equation*}
where  $S:=\inf \{ s\in [0, T]: \frac{p_s}{q_s}(Y_s)>0\}$. This follows from  the pathwise relation 
$$L^r_T\left(\frac{p_{T-\cdot}}{q_{T-\cdot }}(X^Q_{T-\cdot})\right)-L^r_{T-t}\left(\frac{p_{T-\cdot}}{q_{T-\cdot }}(X^Q_{T-\cdot})\right)=L^r_t \left(\frac{p_{\cdot}}{q_{\cdot }}(X^{Q}_.)\right)$$
and the fact that $\left(\frac{p_{T-t}}{q_{T-t}}\left(X^Q_{T-t}\right)\right)_{t\in [0,T]}$ is  a.s. stopped upon hitting $0$,  by Lemma  \ref{DtTexplicit}.
\item[b)]  The limit type assumption  in Corollary  \ref{TV} is not too stringent. Thanks to \eqref{entderiv} and Cauchy-Schwarz inequality, 
 it holds true  for instance if  the matrix $a$ is uniformly  bounded and $H_U(P_0\vert Q_0)<\infty$ for $U(r)=(r-1)^2$, since $ \left| a
\nabla \left[\frac{p_{t}}{q_{t}}\right] \right|=\sup_{|v|\leq 1} (\sigma v)^* (\sigma \nabla  \frac{p_{t}}{q_{t}}) \leq \sqrt{|a| }\sqrt{\nabla^* \left[\frac{p_{t}}{q_{t}} \right]   a
\nabla\left[\frac{p_{t}}{q_{t}}\right]}$.
\end{itemize}
 \end{arem}  
 
 We end this section providing sufficient conditions in order  that $\lim_{t\to\infty}H_U(P_t|Q_t)=0$. The proof of the following result is differed to Appendix  \ref{appvanish}.

\begin{aprop}\label{vanish}  
   Let us assume that the coefficients $\sigma$ and $b$ are time-homogeneous and globally Lipschitz continuous. Then the semigroup associated with the SDE \eqref{diffhom} is Feller. Let us also suppose that \eqref{diffhom} admits an invariant density $p_\infty$, locally Lipschitz and bounded away from $0$ and $+\infty$, and such that  $\int_{\R^d}\frac{p^2_\infty(x)dx}{1+|x|^k}<+\infty$ for some   $k>0$ and that $- \frac{\partial_j
\left(a_{\bullet j}p_{\infty}\right)}{p_{\infty}}$ is the sum of a locally Lipschitz function and a monotone
function. We last suppose that
\begin{equation}
   \forall A>0,\;\exists \varepsilon_A>0,\;\forall |x|\leq A,\;a(x)\geq \varepsilon_A I_d\label{localellip}
\end{equation}
with either $\varepsilon_A$ not depending on $A$ or the second order distribution derivatives $\frac{\partial a^{ij}}{\partial x_i\partial x_j}$ bounded on $\R^d$.
Then,  the tail sigma-field $\cap_{t\geq 0}\sigma(X_r,r\geq t)$ is trivial a.s. w.r.t. the law of $(X^Q_t)_{t\geq 0}$. In particular, if $U(1)=0$, then as soon as $H_U(P_s|Q_s)<+\infty$ for some $s<+\infty$, one has $\lim_{t\to\infty}H_U(P_t|Q_t)=0$.

\end{aprop}

\begin{arem}\label{remtrivsigfield}
The triviality of the tail sigma-field still holds when $(X_t)_{t\geq 0}$  is Feller, has an invariant distribution and a strictly positive transition density $\varphi_t(x,y)$ with respect to the Lebesgue measure which is continuous in $(x,y)$ for each $t>0$ (The continuity implies the strong Feller property, the positivity implies the ergodicity of the invariant measure and combining both, one checks that $(X_t)_{t\geq 0}$ is Harris recurrent. Then one concludes by Theorem 1.3.9 in \cite{Kun'}.) Notice that conditions ensuring the positivity and joint continuity in $(x,y)$ of $\varphi_t(x,y)$ can be found in \cite{friedman} Chapter 9 under uniform ellipticity, and in  \cite{kusustro} Theorem 4.5 under hypoellipticity.
\end{arem}

\section{Dissipation  of the Fisher information and non-intrisic Bakry Emery criterion}\label{secdisfish}
We will from now on focus in the case when  $Q_0(dx)=p_{\infty}(x)dx$ is a stationary probability law for the time-homogeneous Markov diffusion \eqref{diffhom} .
We denote $$I_U(p_s|p_\infty)=\frac{1}{2}\int_{\{\frac{p_{s}}{p_{\infty}}>0\}}U''\left(\frac{p_{s}}{p_{\infty}}\right)\nabla^*\left[
\frac{p_{s}}{p_\infty}\right] a\nabla \left[
\frac{p_{s}}{p_\infty}\right] p_\infty dx$$  
the integral that  appears in the right-hand side of \eqref{entderiv}, and we refer to it as the $U-$ Fisher information. 

Inspired by the famous Bakry-Emery approach, we  want to compute the   derivative of $I_U(p_s|p_\infty)$ with respect to the time variable.  

In all the sequel, we make the following assumptions :
\begin{itemize}
\item[$H4)$]  The drift function $b$ and the matrix $\sigma$ are time-homogeneous and such that $H1)$ holds. Moreover, $b$ (resp. $\sigma$) admits first (resp. second) order derivatives which are locally $\alpha$-H\"older-continuous on $\R^d$ for some $\alpha>0$.
\item[$H5)_{p_\infty}$]  The Markov process defined by  \eqref{diffhom}  has an  invariant density  $p_{\infty}(x)$ and $Q_0(dx)=p_{\infty}(x)dx$. Moreover, $p_{\infty}$  admits derivatives up to the second order which are locally $\alpha$-H\"older-continuous on $\R^d$ for some $\alpha>0$ and $p_{\infty}(x)>0$ for all $x \in \R^d$. 

\item[$H6)_{p_0}^T$]   The initial distribution $P_0$ admits a probability density $p_0$ with respect to the Lebesgue measure. Moreover, we assume that  $H2)_{p_0}$  holds and  that $p_t(x)$ has spatial derivatives up to the second order for each $t>0$, which are continuous in $(t,x) \in (0,T]\times \R^d$ and  bounded and H\"older continuous in $x\in \R^d$ uniformly on $[\delta,T]\times \R^d$ for each $\delta\in (0,T]$. 
\end{itemize}
Let us also introduce some notations :
\begin{itemize}
\item We write   ${\mathbb P}^{T}_{\infty}:= {\mathbb Q}^T$ and $\bar{b}_i:= \bar{b}_i^Q$, $i=1,\dots, d$ . 

% \item  $(A^{-1})_{kl}$ denotes the $(k,l)$
% coordinate of the inverse $A^{-1}$ of an invertible matrix $A$.

\item By possibly enlarging the probability space ${\cal G}_t-{\mathbb P}^{T}_{\infty}$, we introduce a  Brownian motion  $\bar{W}$ such that $Y_t$ solves the stochastic differential equation : 
\begin{equation}\label{diffrev}
dY_t=\sigma(Y_t) d \bar{W}_t+\bar{b}(Y_t)dt,\quad t\in
[0,T]\mbox{ where }\bar{b}_i(y)=-b_i(y)+\frac{\partial_j(a_{ij}(y)p_\infty(y))}{p_\infty(y)}.
\end{equation} 
By assumptions $H4)$  and $H5)_{\infty}$, the coefficients $\sigma$ and $\bar{b}$ are locally Lipschitz so that trajectorial uniqueness holds for this SDE. By the Yamada-Watanabe theorem, one deduces that uniqueness holds for the martingale problem $(MP)_Q$.
\item We write $\rho_t(x):= \frac{p_{T-t}}{p_{\infty}}(x), \  t \in [0,T]$.

\end{itemize}

Notice that $H5)_{p_\infty}$ implies $H2)_Q$ for $Q_0(dx)=p_\infty(x)dx$ and combined with $H4)$, it implies $H3)_Q$. Moreover $H6)_{p_0}^T$ implies $H2)_P$ and $H3)_P$.
Therefore the hypotheses of Theorem \ref{entprod} hold within the present Section. Notice also that,   under $H5)_{\infty}$ and $H6)_{p_0}^T$,  the first order spatial derivatives  of $\frac{p_t}{p_{\infty}}$ are  defined everywhere. Thus,   we may and will assume in the sequel  that   Lemma \ref{DtTexplicit}  b) and   Equation \eqref{entderiv}  hold with the standard gradient $\nabla\frac{p_t}{p_{\infty}}$.  Under $H4)$, if moreover $a$ and $b$ are bounded with $a$ uniformly elliptic, then $H6)_{p_0}^T$ holds for any compactly supported probability density $p_0$, by \cite{friedman} Chapter 9. We refer to \cite{kusustro} for conditions ensuring that $H6)_{p_0}^T$ holds under hypoellipticity.

To compute the dissipation of the $U$-Fischer information, in all the sequel we make   the following regularity assumption on $U$:
\begin{itemize}
\item[$H7)$]  The convex function $U:[0,\infty) \to \R$ is of class $C^4$ on $(0,+\infty)$, continuous on   $[0,+\infty)$ and  satisfies $U(1)=U'(1)=0$.   
\end{itemize}
The assumption that $U'(1)=0$ is inspired in the analysis  on admissible entropies developed in Arnold et al. \cite{AMTU} and  is granted  without modifying the functions $p\mapsto H_U(p\vert p_{\infty})$ and $p\mapsto I_U(p\vert p_{\infty})$ by replacing $U(r)$ by $U(r)-U'(1)(r-1)$ if needed.  Notice that if   $H7)$  holds,  $U(r)$ attains the minimum $0$ at $r=1$ and therefore $U\geq 0$ by convexity.  Following \cite{bakryemery} p.202 (see also \cite{AMTU,chafai}),  we  introduce  an additional assumption on $U$:
\begin{itemize}
\item[$H7')$] $\forall r\in (0,\infty)$, $(U^{(3)}(r))^2\leq \frac{1}{2}U''(r)U^{(4)}(r)$,
\end{itemize}
which is satisfied for instance by $U(r)=r\ln r-(r-1)$ and by $U(r)=(r-1)^2$. Let us recall consequences of $H7)'$ pointed out in \cite{AMTU} (see Remark 2.3 therein) which will be used in proving the following results. \begin{arem}\label{propertiesU}
 Condition $H7')$ implies that  $\left(\frac{1}{U''}\right)''\leq 0$ at points where $U''\not =0$. Since $U''\geq 0$, and excluding the uninteresting case where $U''$  identically  vanishes,  the previous  implies that  $\frac{1}{U''}$ is finite in $[0,\infty)$, and therefore that $U$ is strictly convex. 
 We then deduce from $H7')$ that $U^{(4)}\geq 0$ in $(0,\infty)$. By concavity and positivity of $\frac{1}{U''}$   this function is moreover non decreasing, and we deduce that $U^{(3)}\leq 0$ in $(0,\infty)$. 
\end{arem} 

We do not assume that the entropy function $U$ is $C^4$ on the closed interval $[0,+\infty)$, since we want to deal with $U(r)=r\ln(r)-(r-1)$. That is why we introduce some regularization $U_{\delta}$ indexed by a positive parameter $\delta$ : we chose $U_{\delta}$ such that $U_{\delta}(r)=U(r+\delta)$ for $r\geq 0$ and $U_\delta$ is extended to a $C^4$ function on $\R$. In the next proposition as well as in the remaining of the paper, we will omit the argument $(t,Y_t)$ in order to obtain more compact formulae. 

\begin{aprop}\label{dissipfish}
Under $H4)$, $H5)_{p_\infty}$, $H6)_{p_0}^T$ and $H7)$, one has on the time-interval $[0,T]$
$$ d \left[U''_{\delta} (\rho) \nabla^* \rho a \nabla \rho \right]= tr(\Lambda_{\delta}\Gamma) dt+U''_{\delta}(\rho)\bar{\theta}  dt + d\hat{M}^{(\delta)}\mbox{ with }tr(\Lambda_{\delta}\Gamma)\geq 0\mbox{ under }H7)'$$
and where $\hat{M}^{(\delta)}_t=\int_0^t\partial_k \left[ U''_{\delta}(\rho) \nabla^* \rho a \nabla \rho \right] \sigma_{kr}  d\bar{W}^r_s
$ is a ${\cal G}_t-{\mathbb P}^{T}_{\infty}-$local martingale, % $$d\hat{M}^{(\delta)}:=
% \left\{ 2  \ U_{\delta}''(\rho)  \sigma_{l'i} \  \partial_{l'} \rho   \partial_{k }
% \left[ \partial_l  \rho \sigma_{li}\right]+ \left[\nabla \rho^* a \nabla \rho \right] U^{(3)}_{\delta}(\rho) \partial_k\rho \right\} \sigma_{kr}  d\bar{W}^r =  \partial_k \left[ U''_{\delta}(\rho) \nabla^* \rho a \nabla \rho \right] \sigma_{kr}  d\bar{W}^r, $$ 
\begin{align*}
  \bar{\theta}=
 2 
\bigg\{ &\partial_{l'}\rho \partial_l \rho \left[\frac{1}{4}(\partial_k\sigma_{lj}a_{km}\partial_{m}\sigma_{l'j}-\sigma_{ki}\partial_k\sigma_{lj}\sigma_{mj}\partial_{m}\sigma_{l'i})+\frac{1}{2} \bar{b}_m \partial_m
a_{l l'}+
\frac{1}{2}  \sigma_{l' i} a_{mk}  \partial_{mk}\sigma_{li}
-  a_{m l'}\partial_m
 \bar{b}_l \right] \\&+
  \left[\sigma_{l'i}a_{mk}- \sigma_{ki}a_{ml'} \right] \partial_{l'} \rho
\partial_m\sigma_{li}\partial_{kl} \rho \bigg\},\\   
\end{align*}
and $\Lambda_{\delta}$ and $\Gamma$ are the square matrices defined by
\begin{equation*}
\Lambda_{\delta} :=  \left[\begin{array}{ll} 
   U''_{\delta}(\rho) & U^{(3)}_{\delta}(\rho) \\
 U^{(3)}_{\delta} (\rho) &  \frac{1}{2}U^{(4)}_{\delta} (\rho)\\ \end{array}\right]  \qquad
\Gamma :=  \left[\begin{array}{ll}   \Gamma_{11}&   (\sigma_{\bullet i} \cdot \nabla \rho )  \nabla^*  \rho \  a \nabla(\sigma_{\bullet i} \cdot \nabla \rho)  \\
 (\sigma_{\bullet i} \cdot \nabla \rho )  \nabla^*  \rho \  a \nabla(\sigma_{\bullet i} \cdot \nabla \rho)  & 
  \left|\nabla^*  \rho a \nabla \rho \right|^2 \\ \end{array}\right] 
\end{equation*}
with $\Gamma_{11}=\sum_{i,j=1}^d\left(\sigma_{kj}\sigma_{li}\partial_{kl}\rho+\frac{1}{2}(\sigma_{kj}\partial_{k}\sigma_{li}+\sigma_{ki}\partial_{k}\sigma_{lj})\partial_l\rho\right)^2$
\end{aprop}
The computation of $d \left[U''_{\delta} (\rho) \nabla^* \rho a \nabla \rho \right]$ is postponed to Appendix \ref{pdissipfish}. Let us nevertheless discuss the sign of the term $tr(\Lambda_{\delta}\Gamma)$ which is inspired from \cite{bakryemery} p.202 and also from the term $tr(\mathbf{X}\mathbf{Y})$ in \cite{Arnoldcarlenju} pp 163-164 (see Appendix \ref{compacj} for a detailed comparison with the computations in that paper). Since, by Cauchy Schwarz inequality,
\begin{align*}
   ((\sigma_{\bullet i} \cdot \nabla \rho )  \nabla^*  \rho \  a \nabla(\sigma_{\bullet i} \cdot \nabla \rho))^2&=\left((\sigma^*\nabla\rho)_i(\sigma^*\nabla\rho)_j\left(\sigma_{kj}\sigma_{li}\partial_{kl}\rho+\frac{1}{2}(\sigma_{kj}\partial_{k}\sigma_{li}+\sigma_{ki}\partial_{k}\sigma_{lj})\partial_l\rho\right)\right)^2\\
&\leq \Gamma_{11}\sum_{i,j=1}^d(\sigma^*\nabla\rho)^2_i(\sigma^*\nabla\rho)^2_j=\Gamma_{11}\left|\nabla^*  \rho a \nabla \rho \right|^2,
\end{align*}
the determinant of the matrix $\Gamma$ is nonnegative and this matrix is positive semidefinite. Under $H7')$,  $\Lambda_{\delta}$ is also positive semidefinite and $tr(\Lambda_{\delta}\Gamma)\geq 0$. 
\begin{arem}\label{impnibec}
   In a previous version of this paper, the coefficient $\Gamma_{11}$ was chosen equal to $$\sum_{i,j=1}^d\left(\sigma_{kj}\sigma_{li}\partial_{kl}\rho+\sigma_{kj}\partial_{k}\sigma_{li}\partial_l\rho\right)^2=\sum_{i,j=1}^d(\sigma^*\nabla(\sigma^*\nabla\rho)_i)^2_j=\nabla^*((\sigma\nabla \rho)_i  a   \nabla((\sigma\nabla \rho)_i).$$ We thank Anton Arnold for pointing out to us that the positive semidefiniteness of the matrix $\Gamma$ is preserved under the new choice of this coefficient. Notice that, by symmetry of $\sigma_{kj}\sigma_{li}\partial_{kl}\rho$ in $i$ and $j$,
\begin{align*}
   \sum_{i,j=1}^d\left(\sigma_{kj}\sigma_{li}\partial_{kl}\rho+\sigma_{kj}\partial_{k}\sigma_{li}\partial_l\rho\right)^2-\Gamma_{11}&=\frac{1}{4}\sum_{i,j=1}^d\left((\sigma_{kj}\partial_{k}\sigma_{li}-\sigma_{ki}\partial_{k}\sigma_{lj})\partial_l\rho\right)^2\\&=\frac{1}{2}(\partial_k\sigma_{lj}a_{km}\partial_{m}\sigma_{l'j}-\sigma_{ki}\partial_k\sigma_{lj}\sigma_{mj}\partial_{m}\sigma_{l'i})\partial_l\rho\partial_{l'}\rho
\end{align*}
is a nonnegative quadratic form applied to $\nabla\rho$ which implies that the Bakry Emery criterion below improves upon the one of  the previous version.
\end{arem}
% Notice that one could preserve the positive semidefiniteness of the matrix $\Gamma$ when replacing $\Gamma_{11}$ by the smaller coefficient $\sum_{i=1}^{d'}\left(\nabla^*  \rho \  a \nabla(\sigma_{\bullet i} \cdot \nabla \rho)\right)^2/\left|\nabla^*  \rho a \nabla \rho \right|$, which amounts to replace the squared norms of the vectors $\sigma^*\nabla(\sigma_{\bullet i} \cdot \nabla \rho)$  by the ones of their orthogonal projection on $\sigma^*\nabla\rho$. Unfortunately, we have not been able to take advantage of this possibility.

We introduce one  last  assumption on the density flow $\rho_t=\frac{p_{T-t}}{p_\infty}$ :

$H6')_{p_0}^T$   For each $T'\in (0,T)$ the following integrals are finite:
\begin{itemize}
\item $\int_0 ^{T'}\left| U^{(3)}(\rho) \vee-1 \right|^2  | \nabla^*\rho  a\nabla\rho  |^3 p_{\infty}(x)dx dt $
\item    $ \int_0 ^{T'} \left( U''(\rho) \wedge 1\right)^2   \nabla^*( \nabla^*\rho a\nabla\rho  ) a \nabla ( \nabla^*\rho a\nabla\rho )   p_{\infty}(x)dx dt$
\item $  \int_0 ^{T'} (U''(\rho) \wedge 1 ) \big[  \left| (\sigma_{l'i}a_{m\bullet}- \sigma_{\bullet i}a_{ml'}  ) \partial_m\sigma_{li} \right| +   \left| \partial_k \left( \left[\sigma_{l'i}a_{mk}- \sigma_{ki}a_{ml'} \right] \partial_m\sigma_{li} \right)\right|
 \big] | \partial_{l'} \rho | | \partial_{l} \rho | p_{\infty}(x)dx dt  $ 
  \item $  \int_0 ^{T'} (U''(\rho) \wedge 1 ) \big[  \left| (\sigma_{l'i}a_{mk}- \sigma_{k i}a_{ml'}  ) \partial_m\sigma_{li} ( \partial_l \rho \partial_ k \ln p_{\infty}  +  \partial_{lk} \rho) \right| 
 \big] | \partial_{l'} \rho | p_{\infty}(x)dx dt  $   
\end{itemize}
We also denote by $H6)_{p_0}^{\infty}$ (resp. $H6')_{p_0}^{\infty}$) the assumption that $H6)_{p_0}^T$ (resp. $H6')_{p_0}^T$) holds for each $T>0$.

%  Notice that in the case that $a\geq c I_d$ for some $c>0$, the third integral converges  always as it can be upper bounded by the Fisher information. 

\begin{thm}\label{thdissipentro}
Let $\Theta$ denote the $d\times d$ symmetric matrix defined by 
\begin{equation*}
\begin{split}
   \Theta_{ll'}&=- \frac{1}{2} b_m \partial_m a_{l l' } + \frac{1}{2} ( a_{kl'}\partial_k b_l+  a_{kl}\partial_k b_{l'} )  -\frac{1}{4} a_{mk}\partial_{mk} a_{l l' }  -\frac{1}{2} (a_{kl'}\partial_{kj}a_{lj}+a_{kl}\partial_{kj}a_{l'j} ) \\ 
 & - a_{kl}a_{jl'}\partial_{kj}\ln (p_{\infty}) -\frac{1}{2} (a_{kl} \partial_k a_{l'j}  + a_{kl'} \partial_k a_{lj} )\partial_j   \ln (p_{\infty}) - \frac{1}{4} ( a_{mk} \partial_m \sigma_{li}\partial_k \sigma_{l'i}+\sigma_{ki}\partial_k\sigma_{lj}\sigma_{mj}\partial_m\sigma_{l'i}) \\
 & +  \frac{1}{2}  \sigma_{ki} (\partial_m\sigma_{li}a_{ml'}+  \partial_m\sigma_{l' i} a_{ml})\partial_k\ln(p_\infty) +  \frac{1}{2}   \partial_k[\sigma_{ki} (\partial_m\sigma_{li}a_{ml'}+  \partial_m\sigma_{l' i} a_{ml})  ] 
 \end{split}
\end{equation*}
and assume  that $\Theta(x)$ is $p_{\infty}(x)dx-  a.e.$ positive  semidefinite. Then, under $H4)$, $H5)_{p_\infty}$, $H6)_{p_0}^T$ $H6')_{p_0}^T$, $H7)$ and $H7')$, for a.e. $t\in [0,T]$ one has 
\begin{align}
   \frac{d}{dt}\int_{\rho_t>0} U''(\rho_t)[\nabla^*\rho_t a\nabla\rho_t ] p_\infty dx\geq  2 \int_{\rho_t>0}U''(\rho_t)\nabla^*\rho_t\Theta\nabla\rho_t p_\infty dx\label{minodissip2}.
\end{align}
If moreover $I_U(p_0|p_\infty)<+\infty$, $H6)_{p_0}^{\infty}$ and $H6')_{p_0}^{\infty}$ hold and the matrix $\Theta$ satisfies the non-intrinsic Bakry-Emery criterion
\begin{itemize}
\item[NIBEC)] $\exists \lambda>0,\;\forall x\in \R^d,\;\Theta(x)\geq \lambda a(x),$
\end{itemize}
then 
$\forall t\geq 0$, $I_U(p_t|p_\infty)\leq e^{-2\lambda t}I_U(p_0|p_\infty)$ and the non-increasing function $t\mapsto H_U(p_t|p_\infty)$ also converges at exponential rate $2\lambda$ to its limit as $t\to\infty$.
\end{thm}
\begin{arem}\label{remcassym}\begin{itemize}
   \item  The matrix $\Theta$ and therefore our Bakry-Emery criterion are non-intrinsic in  the sense that they cannot in general be written in terms of the diffusion matrix $a$ only, without making explicit use of $\sigma$. This is because we have got rid of the nonnegative term $tr(\Lambda_\delta\Gamma)$ which appears in the first equation in Proposition \ref{dissipfish} and involves the non-intrisic term $\Gamma_{11}$.
\item In case $a=2\nu I_d$ and $b=-(\nabla V+F)$ with $F$ such that $\nabla.(e^{-V/\nu}F)=0$, then $p_\infty\propto e^{-V/\nu}$, $\bar{b}=-b+2\nu\nabla\ln p_\infty=-\nabla V+F$ and $\Theta=\nu(2\nabla^2 V-\nabla F-\nabla F^*)$. For the canonical choice  $\sigma=\sqrt{2\nu}I_d$, condition  {\it NIBEC)}    therefore   writes    
$\exists \lambda >0,\;\forall x\in\R^d,\;\nabla^2 V(x)-\frac{\nabla F+\nabla F^*}{2}(x)\geq \lambda I_d$ which is exactly condition (A2) in the introduction of \cite{Arnoldcarlenju}, page 158.  

\end{itemize}

\end{arem}
The proof of \eqref{minodissip2} is postponed to Appendix \ref{pthdissipentro}. Let us deduce the last assertion of Theorem \ref{thdissipentro}. Reverting time in \eqref{minodissip2} and using {\it NIBEC)}, one obtains that for $r\geq 0$,
$$\frac{d}{dr}I_U(p_r|p_\infty)\leq -2\lambda I_U(p_r|p_\infty).$$
Hence $\forall r\geq t\geq 0,\;I_U(p_r|p_\infty)\leq e^{-2\lambda (r-t)} I_U(p_t|p_\infty)$. Since by Theorem \ref{entprod}, one has $\frac{d}{dr}H_U(p_r|p_\infty)=-I_U(p_r|p_\infty)$, we deduce  that 
\begin{equation}
   0\leq H_U(p_t|p_\infty)- \lim_{r\to\infty} H_U(p_r|p_\infty)= \int_t^\infty I_U(p_r|p_\infty)dt\leq \frac{I_U(p_t|p_\infty)}{2\lambda}\leq \frac{e^{-2\lambda t}I_U(p_0|p_\infty)}{2\lambda}.\label{comphi}
\end{equation}

We   deduce

\begin{thm}\label{thsl}Assume $H4)$, $H5)_{p_\infty}$, $H6)_{p_0}^{\infty}$ $H6')_{p_0}^{\infty}$, $H7)$ and $H7')$, that the matrix $\Theta(x)$ is $p_{\infty}(x)dx-  a.e.$ positive  semidefinite, that the diffusion matrix $a$ is locally uniformly strictly positive definite and that $H_U(p_s|p_\infty)<+\infty$ for some $s\geq 0$. Then $H_U(p_t|p_\infty)$ converges to $0$ as $t\to\infty$. 
Moreover, under $NIBEC)$, for $t>s$, one has the convex Sobolev inequality 
\begin{align}
 &H_U(p_t|p_\infty)\leq \frac{1}{2\lambda}I_U(p_t|p_\infty), \label{convsobineq}\\
\mbox{and   }&\forall t\geq s,\;H_U(p_t|p_\infty)\leq e^{-2\lambda (t-s)}H_U(p_s|p_\infty).\label{decexpo}
\end{align}
\end{thm}

\begin{adem}Reverting time in \eqref{minodissip2}, we obtain that $t\mapsto I_U(p_t|p_\infty)$ is non-increasing. When $H_U(p_s|p_\infty)$ is finite for some $s\geq 0$, writing \eqref{entderiv} on the interval $[s,T]$ in place of $[0,T]$ with arbitrarily large $T$, we deduce that $I_U(p_t|p_\infty)$ is finite on $(s,+\infty)$ and tends to $0$ as $t\to\infty$. When $a$ is locally uniformly strictly positive definite, the beginning of the proof of Theorem 2.5  \cite{Arnoldcarlenju} (before Part(a)), ensures that $p_t$ tends to $p_\infty$ in $L^1(\R^d)$. As a consequence, in the notations of the introduction, $\E\left|\frac{dP_t}{dQ_t}(X^Q_t)-1\right|$ tends to $0$ as $t\to\infty$ and therefore the a.s. limit $\E\left(\frac{dP_t}{dQ_t}(X^{Q}_t)\bigg|\cap_{s\geq 0}{\cal F}_s\right)$ of $\frac{dP_t}{dQ_t}(X^Q_t)$ is equal to $1$. By \eqref{limentrel}, one concludes that $H_U(p_t|p_\infty)$ tends to $U(1)=0$.\\
Under $NIBEC)$, for $t>s$, $I_U(p_t|p_\infty)<+\infty$ and reasoning like in the derivation of 
\eqref{comphi}, one obtains \eqref{convsobineq}. This implies that
$$\frac{d}{dt}H_U(p_t|p_\infty)=-I_U(p_t|p_\infty)\leq -2\lambda H_U(p_t|p_\infty)$$
from which the last assertion follows readily.
\end{adem}
\begin{arem}
In view of \eqref{limentrel} and Remark \ref{remtrivsigfield}, the local uniform strict positive definiteness assumption on the diffusion matrix $a$ may be replaced by some hypoellipticity assumption, in order to ensure that $H_U(p_t|p_\infty)$ tends to $0$ as $t\to\infty$ at exponential rate $2\lambda$ as soon as $H_U(p_s|p_\infty)<\infty$  for some $s\geq 0$. By the last step of the proof of Theorem \ref{thsl}, this implies \eqref{convsobineq} and \eqref{decexpo} under $NIBEC)$.
\end{arem}

\section{Examples}\label{EXs}

%\subsection{}  
Consider the reversible diffusion process in $\R^2$ with coefficients given for each $(x_1,x_2)\in \R^2$ by 
 $$a(x_1,x_2)=I_2, \quad\mbox{ and } \quad  b(x_1,x_2)=-\nabla V(x_1,x_2),$$
 where $V$ is the globally $C^2$ convex potential 
\begin{equation*}
   V(x_1,x_2):=|x_1|^2+|x_1-x_2|^{2+\alpha} + |x_2|^{2+\alpha}% \mbox{ with }v(z)=\begin{cases}z^{2+\alpha}\mbox{ if }z\in[0,1]\\
%1+(2+\alpha)(z-1)\frac{(1+\alpha)z+(1-\alpha)}{2}\mbox{ if }z\geq 1\end{cases}.
\end{equation*}
 for some $\alpha \in (0,1)$.  
The invariant measure is $p_\infty\propto e^{-2V}$, and we
 have
\begin{equation*}
\begin{split}
\partial_1 V =& 2x_1 +(2+\alpha)sign(x_1-x_2)|x_1-x_2|^{1+\alpha} % \1_{\{|x_1-x_2|\leq 1\}}+((1+\alpha)|x_1-x_2|-\alpha)\1_{\{|x_1-x_2|>1\}}] 
 \\
\partial_2 V = & (2+\alpha)sign(x_2)|x_2|^{1+\alpha}%\1_{\{|x_2|\leq 1\}}+((1+\alpha)|x_2|-\alpha)\1_{\{|x_2|>1\}}]
+ (2+\alpha)sign(x_2-x_1)|x_2-x_1|^{1+\alpha}%\1_{\{|x_2-x_1|\leq 1\}}+((1+\alpha)|x_2-x_1|-\alpha)\1_{\{|x_2-x_1|>1\}}]
\\
\end{split}
\end{equation*}
and
\begin{equation*}
\nabla^2 V = \left( \begin{array}{cc}
 2 &    0  \\
0 & (2+\alpha)(1+\alpha)|x_2|^{\alpha} \\
\end{array} \right) +  (2+\alpha)(1+\alpha) |x_1- x_2|^{\alpha} \left( \begin{array}{cc}
 1 &    -1 \\
-1  & 1  \\
\end{array} \right). 
\end{equation*}
The drift  $b=-\nabla V$ is locally Lipschitz continuous. Moreover, $(x_1,x_2).\nabla V(x_1,x_2)\geq 0$ and $\sqrt{\partial_{ik}V\partial_{ik}V(x_1,x_2)}\leq C\sqrt{1+|x_2|^{2\alpha}+|x_1-x_2|^{2\alpha}}$ so that $\limsup_{|(x_1,x_2)|\to+\infty}\frac{\sqrt{\partial_{ik}V\partial_{ik}V(x_1,x_2)}}{V(x_1,x_2)}=0$. Last $\Delta V(x_1,x_2)\leq C(1+|x_2|^{\alpha}+|x_1-x_2|^{\alpha})$ whereas\begin{align*}
   |\nabla V|^2(x_1,x_2)\geq &(2|x_2|+(2+\alpha)|x_1-x_2|^{1+\alpha})^2{\mathbf 1}_{sign(x_2)\neq sign(x_2-x_1)}\\&+(2+\alpha)^2(|x_2|^{1+\alpha}+|x_1-x_2|^{1+\alpha})^2{\mathbf 1}_{sign(x_2)=sign(x_2-x_1)}
\end{align*}since $sign(x_2)\neq sign(x_2-x_1)$ iff $x_1\geq x_2\geq 0$ or $x_1\leq x_2\leq 0$. Therefore $\limsup_{|(x_1,x_2)|\to+\infty}\frac{\Delta V}{|\nabla V|^2}(x_1,x_2)=0$ and, by Remark \ref{condthmprincip} b), $H1)''$ is satisfied.

 The classic Bakry-Emery criterion fails since $\nabla^2V(0,0)$ is singular but a logarithmic Sobolev inequality can be obtained by the perturbative argument of Holley and Stroock \cite{HS}.   The potential $V$ is also a particular case of the examples considered by Arnold, Carlen and Ju in the Section 3 of \cite{Arnoldcarlenju}.  We notice that in order to check that $p_\infty$ satisfies the convex Sobolev inequality \eqref{convsobineq}, they first modify the Fokker-Planck equation by adding a non-symmetric drift term $F$ as described 
 in Remark \ref{remcassym} ii) above. Exponential convergence to $0$ of $H_U(p_t|p_\infty)$ for the solution $p_t$ of the original Fokker-Planck equation is only deduced in a second step. 
 
 It is nevertheless of interest to see how  our non-intrisic Bakry Emery criterion allows us to prove directly  that $p_\infty$  satisfies the convex Sobolev inequality \eqref{convsobineq} and that $H_U(p_t|p_\infty)$ converges exponentially to $0$. In contrast to  \cite{Arnoldcarlenju} we  modify the stochastic differential equation associated with the diffusion processes, by  changing the square root $\sigma$  of the diffusion matrix, but not the law of its solution or the associated Foker-Planck equation. We consider   
$$\sigma=\left( \begin{array}{cc}
 \cos \phi  &   \sin \phi   \\
-    \sin \phi  &  \cos \phi    \\
\end{array} \right)$$
for a function $\phi:\R^2\to \R^2$ of class $C^2$ to be chosen later. We obtain after some  computations  
\begin{align*}
 \Theta =\nabla^2 V  -\frac{1}{4} |\nabla \phi |^2 I_2&-\frac{1}{4} \left( \begin{array}{cc}
 (\partial_{2}\phi)^2  &   -\partial_{1}\phi\partial_{2}\phi      \\
   -\partial_{1}\phi\partial_{2}\phi     &  (\partial_{1}\phi)^2 \\
\end{array} \right)+  \left( \begin{array}{cc}
 \partial_{12}\phi  &   \frac{\partial_{22}\phi -\partial_{11}\phi}{2}      \\
   \frac{\partial_{22}\phi-\partial_{11}\phi} {2}     &  -\partial_{12}\phi     \\
\end{array} \right)   \\ &+ \left( \begin{array}{cc}
 -2\partial_1{\phi}\partial_2 V   &   \partial_1{\phi}\partial_1 V -\partial_2{\phi}\partial_2 V     \\
 \partial_1{\phi}\partial_1 V -\partial_2{\phi}\partial_2 V  &  2\partial_2{\phi}\partial_1 V     \\
\end{array} \right) 
\end{align*}

We now consider  a parameter $\varepsilon>0$ which will be chosen small and a $C^2$ function $\varphi:\R\to \R$ such that $\varphi(s)=s$ if $|s|\leq 1$ and $\varphi(s)=0$ if $|s|\geq 2$. Then, we define
$$\phi(x_1,x_2)=-\varepsilon \varphi_{\varepsilon}(x_1) \varphi_{\varepsilon}(x_2), \quad (x_1,x_2)\in \R^2$$
where $\varphi_{\varepsilon}(s)= \varepsilon \varphi(s/{\varepsilon})$. Notice that 
$$\varphi_{\varepsilon} =O(\varepsilon), \quad \varphi_{\varepsilon}'' =O(1/\varepsilon), \quad\mbox{ and }\,   \varphi_{\varepsilon}' =\left\{ \begin{array}{ll} 1 &\mbox{if } |s|\leq \varepsilon, \\ O(1)&\mbox{if }\varepsilon<|s|<2\varepsilon,\\ 0  &\mbox{if }  |s|\geq 2\varepsilon .\\ \end{array}\right. $$
 Then,  defining $B_{\varepsilon}:=\{(x_1,x_2)\in \R^2 \, s.t.  \, |x_1|\vee|x_2| \leq \varepsilon \}$  and $C_{\varepsilon}:= B_{2\varepsilon}\backslash B_{\varepsilon}$, we have 
$$\partial_1\phi (x_1,x_2), \partial_2\phi (x_1,x_2)=\left\{ \begin{array}{ll}  O(\varepsilon^2)  &\mbox{if } (x_1,x_2)\in B_{2\varepsilon},  \\
0  &\mbox{if } (x_1,x_2)\in B_{2\varepsilon}^c ,  \\ \end{array} \right. $$  $$ \partial_{12}\phi (x_1,x_2)=\left\{ \begin{array}{ll}  -\varepsilon  &\mbox{if } (x_1,x_2)\in B_{\varepsilon},  \\
O(\varepsilon) &\mbox{if } (x_1,x_2)\in C_{\varepsilon} , \\ 
0 &\mbox{if }  (x_1,x_2)\in B_{2\varepsilon}^c ,  \end{array} \right.  $$
$$  \frac{1}{2}( \partial_{11}\phi (x_1,x_2)-  \partial_{22}\phi (x_1,x_2) )=\left\{ \begin{array}{ll}  0  &\mbox{if } (x_1,x_2)\in B_{\varepsilon},  \\
O(\varepsilon) &\mbox{if } (x_1,x_2)\in C_{\varepsilon}, \\ 
0 &\mbox{if }  (x_1,x_2)\in B_{2\varepsilon}^c, \end{array} \right.  $$
and   $\partial_1 V =O(\varepsilon), \partial_2 V =O(\varepsilon^{1+\alpha})$  on $B_{2\varepsilon}$.    It follows that
$$\Theta =\nabla^2 V +  \left( \begin{array}{cc}
-\varepsilon  &   0   \\
 0    &   \varepsilon      \\
\end{array} \right)   + O(\varepsilon^3 ) \geq   \left( \begin{array}{cc}
2 -\varepsilon  &   0   \\
 0    &   \varepsilon      \\
\end{array} \right)  + O(\varepsilon^3 ) \quad \mbox{ on }B_{\varepsilon}.$$

Next, the smallest eigenvalue of $\nabla^2 V(x_1,x_2)$,  is given by
\begin{equation}
   \gamma_{-}:= 1+  \kappa_1   +\kappa_2/ 2 - \sqrt{ 
 1+  \kappa_1^2 -  \kappa_2 + \kappa_2^2/4  }\geq 1+{\kappa_2}/{2}-\sqrt{({\kappa_2}/{2}-1)^2}=\kappa_2\wedge 2\label{minogamma-}
\end{equation}
with $\kappa_1=\kappa_1(x_1,x_2):= (2+\alpha)(1+\alpha)|x_1- x_2|^{\alpha} $ and $\kappa_2 =\kappa_2(x_1,x_2):=  (2+\alpha)(1+\alpha)|x_2|^{\alpha} $.
Since $\gamma_{-}=\kappa_1   +\kappa_2+O(\kappa_1^2+\kappa_2^2)$ as $\kappa_1^2+\kappa_2^2\to 0$ and $|x_2|^{\alpha}+|x_1-x_2|^{\alpha}\geq (|x_2|+|x_1-x_2|)^\alpha\geq |x_1|^\alpha$, we deduce that on   $ C_{\varepsilon}$, 
$$\Theta=\nabla^2 V + O(\varepsilon)\geq (2+\alpha)(1+\alpha)\varepsilon^{\alpha} I_2 +o(\varepsilon^\alpha).$$
Last, % as $\kappa_1^2+\kappa_2^2\to+\infty$, $\gamma_{-}=1+  \kappa_1   +\kappa_2/ 2 -\sqrt{\kappa_1^2+\kappa_2^2/4}+\frac{\kappa_2}{\sqrt{4\kappa_1^2+\kappa_2^2}}+o(1)$ which ensures that $\liminf_{\kappa_1^2+\kappa_2^2\to+\infty}\gamma_-=1$. Since $\lim_{x_1^2+x_2^2\to+\infty}\kappa_1^2+\kappa_2^2=+\infty$ and
by \eqref{minogamma-}, $\inf_{(x_1,x_2)\in B_{2\varepsilon}^c}\gamma_-\geq ((2+\alpha)(1+\alpha)(2\varepsilon)^\alpha)\wedge 2>0$. We conclude that for $\varepsilon$ small enough {\it NIBEC)} holds.

%\subsection{}

\bigskip

We next study a related second example of application of our criterion, 
where $\nabla^2V$ is singular on a ball with positive radius.
Once again, the perturbative argument of Holley Stroock [12] also ensures that
a logarithmic Sobolev inequality holds for this choice of potential.

 Let $v$ be a convex $C^2$ function which vanishes on $[-\frac{1}{4},\frac{1}{4}]$ and such that $v''=2$ on $(-\infty,\frac{1}{2}]\cup[\frac{1}{2},+\infty)$. We set $v_\varepsilon(s)=\varepsilon^2v\left(\frac{s}{\varepsilon}\right)$ and $V_\varepsilon(x_1,x_2)=x_1^2+v_\varepsilon(x_2)+v_\varepsilon(x_1-x_2)$.
For $\varepsilon<\frac{1}{3}$, let $\varphi_\varepsilon$ be a $C^2$ function such that
\begin{equation*}
   \varphi_\varepsilon(s)=\begin{cases}s\mbox{ when }|s|\leq \varepsilon\\0\mbox{ when }|s|\geq 1\end{cases}
 \end{equation*} 
and such that $\frac{-2\varepsilon}{1-\varepsilon}\leq \varphi_\varepsilon'\leq 1$, $|\varphi_\varepsilon|\leq 2\varepsilon$ and $|\varphi''_\varepsilon|\leq C$ where $C$ is a constant not depending on $\varepsilon$.
We set $\phi(x_1,x_2)=-\varphi_\varepsilon(x_1)\varphi_\varepsilon(x_2)$ so that $-1\leq \partial_{12}\phi(x_1,x_2)\leq \frac{2\varepsilon}{1-\varepsilon}$ with the first inequality being an equality on $B_\varepsilon$. We have  $|\partial_{22}\phi-\partial_{11}\phi|\leq 4C\varepsilon$ and $|\nabla \phi|=O(\varepsilon)$. As a consequence, $\Theta=\hat{\Theta}+O(\varepsilon)$ where
$$\hat{\Theta}=\left( \begin{array}{cc}
2+v_\varepsilon''(x_1-x_2)+\partial_{12}\phi(x_1,x_2) &   -v_\varepsilon''(x_1-x_2)    \\
-v_\varepsilon''(x_1-x_2)  &  v_\varepsilon''(x_2)+v_\varepsilon''(x_1-x_2)-\partial_{12}\phi(x_1,x_2)      
\end{array} \right).$$
On $B_\varepsilon$, we have $\partial_{12}\phi(x_1,x_2)=-1$ and $\hat{\Theta}\geq I_2$. If  $|x_2|\geq \frac{\varepsilon}{2}$, then $v''_\varepsilon(x_2)=2$ so that $\hat{\Theta}\geq (2-1)\wedge\left(2-\frac{2\varepsilon}{1-\varepsilon}\right)I_2$. When $|x_2|\leq \frac{\varepsilon}{2}$ and $|x_1|>\varepsilon$,   $|x_1-x_2|\geq \frac{\varepsilon}{2}$ holds so that $v''_\varepsilon(x_1-x_2)=2$ and
$$\hat{\Theta}\geq \left( \begin{array}{cc}
4+\partial_{12}\phi &   -2   \\
-2 &  2-\partial_{12}\phi  
\end{array} \right)\geq \left(3-\sqrt{5+2\partial_{12}\phi-(\partial_{12}\phi)^2}\right)I_2\geq \left(3-\sqrt{5+\frac{4\varepsilon}{1-\varepsilon}}\right)I_2.$$
We conclude that
$$\forall \lambda\in (0,3-\sqrt{5}),\mbox{ for $\varepsilon>0$ and small enough },\forall x\in\R^d,\;\Theta(x)\geq \lambda I_2.$$

\appendix

\section{Proofs of the main results of Section \ref{secentdiss}}\label{proofssecentdiss}
 \subsection{Proof of Lemma \ref{DtTexplicit}}\label{pDtTexplicit}
 
 The proof  of part a) relies on the following technical result:

  \begin{alem}\label{gradlog}
Assume  that $H1)$ , $H2)_P$ and  $H3)_P$ hold.
\begin{itemize}
\item[i)] For each $i=1\dots,d$ and a.e. $t\in (0,T]$,  the distribution $[a_{ij}(t,\cdot)\partial_j p_{t}]:=\partial_j (a_{ij}(t,\cdot) p_{t})- p_t\partial_j a_{ij}(t,\cdot)  $ is a function in $L^1_{loc}(dx)$ and,  as a Radon  measure in $[0,T]\times \R^d$, one has $[a_{ij}(t,\cdot)\partial_j p_{t}](x)dx\ dt \ll p_t(x)dx\ dt$. A measurable in $(t,x)$ version of the Radon-Nikodyn density  is given by  $[a_{i j}(t,\cdot ) \partial_j p_t] (x) / p_t(x)$.  Moreover,
   there exists a measurable function $(t,x)\mapsto K^p(t,x)\in \R^d$ such that for each $i=1\dots,d$
   $$[a_{i j}(t,\cdot) \partial_j p_t] (x) / p_t(x)=  a_{i \bullet}(t,x) K^p(t,x) , \ p_t(x)dx \ dt \ a.e.$$
where $a_{i\bullet}$ denotes the row vector $(a_{i1},\hdots,a_{id})$.
   
\item[ii)] If moreover $H2)_Q$,  $H3)_Q$ and  $P_0\ll Q_0$  hold,
 one has $[a_{ij}(t,\cdot)\partial_j p_{t}](x)dx\ dt \ll q_t(x)dx \ dt$ and  $[a_{i j}(t,\cdot) \partial_j p_t] (x) / q_t(x)$  is a  measurable in $(t,x)$ version of the Radon-Nikodyn derivative. Furthermore, it  holds $p_{T-t}(x)dx \ dt$  (but not necessarily $q_{T-t}(x)dx \ dt$) a.e. that 
\begin{equation*}
\begin{split}
\bar{b}_i^P(t,x)-\bar{b}_i^Q(t,x)= & [\bar{a}_{ij}(t,\cdot )\partial_j p_{T-t}](x) /  p_{T-t}(x) -[\bar{a}_{ij}(t,\cdot)\partial_j q_{T-t}](x)/ q_{T-t}(x) \\ 
= & \bar{a}_{i \bullet}(t,x)( K^p(T-t,x)- K^q(T-t,x)), \\
\end{split}
\end{equation*} 
and $q_{T-t}(x)dx \ dt$ (and thus $p_{T-t}(x)dx \ dt$) a.e.  that 
\begin{equation*}
\begin{split}
\frac{p_{T-t}(x)}{q_{T-t}(x)}(\bar{b}_i^P(t,x)-\bar{b}_i^Q(t,x))= & \frac{p_{T-t}(x)}{q_{T-t}(x)}\bar{a}_{i \bullet}(t,x)( K^p(T-t,x)- K^q(T-t,x)). \\
\end{split}
\end{equation*}  
\end{itemize}
\end{alem}

\begin{adem}  The   Lipschitz character of $a$ (following from $H1)$) ensures that  $a$ has  a.e. defined  spatial derivatives of order $1$ in $L^{\infty}_{loc}([0,T]\times\R^d)$. Thus,  the distribution $a_{i j} (t,\cdot)\partial_j p_t$ is a function
 in $L^1_{loc}([0,T]\times\R^d)$ under  $H3)_P$. This implies, by Lemma A.2 in \cite{MNS}  (see also Lemma A.2 in \cite{HP}),  that   $a_{i j}(t,x) \partial_j p_t (x)$ vanishes   a.e. on $\{x: p_t(x)=0\}$. This fact easily yields the remaining assertions, except the existence of the functions $K^p$ or $K^q$, which we establish in what follows. 
 
 We will on one hand use the fact asserted in the proof of Lemma A.2 in \cite{MNS} that, for a.e. $t>0$ and  each bounded open set $O$, $a_{i j}(t,x) \partial_j p_t (x)$ is the $\sigma(L^1(O),L^{\infty}(O))$-weak limit of some subsequence of  $a_{i j}(t,x) \partial_j [\rho_n*p_t] (x)$, for compactly supported regularizing kernels   $\rho_n(x)=n^d\rho(nx)$. It is indeed shown in  Lemma A.1 in \cite{HP}  that   for a suitable bounded sequence $\alpha_n>0$, $\alpha_n^{-1} |x|\ | \nabla \rho_n(x)|$   is again a regularizing kernel. The local Lipschitz character of $a$ then yields  the domination   $\forall x\in O$, $|a_{i j}(t,x) \partial_j [\rho_n*p_t] (x)|\leq |\rho_n* \partial_j (a_{ij}(t,\cdot) p_{t}) (x)|+ C \alpha_n^{-1}\int |x-y|\ | \nabla \rho_n(x-y)|p_t(y)dy  $,  the right-hand side being,  by the previous,  an    $L^1(O)$-converging sequence.  Weak compactness is then provided by the Dunford-Pettis criterion, and the limit is identified integrating by parts against smooth test functions
 compactly supported  in $O$.   On the other hand, diagonalizing the symmetric positive semidefinite matrix $(a_{ij}(t,x))=[u_1(t,x),\dots,u_d(t,x)]  \Lambda(t,x) [u_1(t,x),\dots,u_d(t,x)] ^*$ provides  orthonormal vectors $(u_i(t,x))_{i=1}^d$ and the corresponding eigenvalues and diagonal components  $(\lambda_i(t,x))_{i=1}^d$ of $\Lambda(t,x)$,  that are measurable as  functions  of   $(t,x)$. 

We  take   $O$ as before and   $a_{i j}(t,x) \partial_j [\rho_n*p_t] (x)$ to be the subsequence described above. Defining the vectorial functions $w^{(n)}:= [u_1,\dots,u_d] ^*\nabla [\rho_n*p_t]$ and $v_k=sign(u_k^*[a \nabla p] ) u_k, \  k=1,\dots,d$, we  have
\begin{equation*}
\int_{O\cap \{\lambda_k=0\} } | v_k^* [ a \nabla p_t] | =\lim_{n\to \infty} 
\int_{O\cap \{\lambda_k=0\} } v_k^* [ a \nabla  [\rho_n*p_t]  ]  =  \lim_{n\to \infty} 
\int_{O\cap \{\lambda_k=0\} }  \lambda_k w^{(n)}_k sign(u_k^*[a \nabla p_t] )=0,
\end{equation*}
since $a \nabla  [\rho_n*p_t]=\sum_{j=1}^d \lambda_j w^{(n)}_j u_j $ by the spectral decomposition of $a$.  Consequently, for each $t$ and  a.e. $x\in \R^d$, the vector $[a(t,x) \nabla p_t(x)]$  belongs to the linear space  $	\big\langle( u_i(t,x))_{i=1,\dots, d ; \lambda_i(t,x)\not = 0}\big\rangle$.  Denote now by $w=(w_j)_{j=1}^d:=( u_j^*a\nabla p_t)_{j=1}^d $  the coordinates of $a\nabla p_t$ w.r.t. the orthogonal basis $( u_j(t,x))_{j=1,\dots, d}$, so that  $w$  is a measurable function of $(t,x)$.   If we moreover denote by  $\overline{\Lambda}$ the diagonal matrix with diagonal coefficients $\lambda^{-1}_j\mathbf{1}_{\lambda_j\not =0}, j=1,\dots ,d$, and set  $v:= [u_1,\dots,u_d] \overline{\Lambda}w$, then
$$ a v= [u_1,\dots,u_d] \Lambda [u_1,\dots,u_d]^* [u_1,\dots,u_d] \overline{\Lambda}w = [u_1,\dots,u_d] \Lambda\overline{\Lambda}w= [u_1,\dots,u_d] w $$ 
since  $w=(w_j\mathbf{1}_{\lambda_j\not =0})_{j=1}^d$.  That is, $(t,x)\mapsto v(t,x)\in \R^d$ is a measurable function such that for almost every $t\in [0,T] $ and each $i$, $a_{i\bullet}(t,x)v(t,x)=[a_{ij}\partial_j p_t(x)], \ dx$ a.e.  Finally,  $K^p(t,x):=v(t,x)/p_t(x)\mathbf{1}_{p_t(x)>0}$ has the required properties.  \end{adem}

\medskip

We can now take  $\nabla \ln \frac{p_t}{q_t}(x)$  to be an arbitrary representant of the  equivalence class of the function $ K^p(t,x)- K^q(t,x)$ under the relation $f(t,x)-g(t,x)\in Ker(a(t,x)), \ p_t(x)dx\ dt\  a.e.$
The  identity in Lemma \ref{DtTexplicit} a) is then satisfied by construction.

 The proof of part b) of Lemma  \ref{DtTexplicit}  firstly relies  on a martingale representation property ensured by the extremality assumption, according to Theorem 12.21 in \cite{Jac}:

\begin{alem}\label{reprmart} Assume that $H1), H2)_Q$ and $H3)_Q$ hold.  
For each $i=1,\dots,d$, 
$$M_t^i:=Y_t^i-Y_0^i-\int_0^{t} \bar{b}_i^Q(s,Y_s)ds, \ t\in [0,T]$$
 is a continuous local
martingale with respect to  $\Q^T$ and $({\cal G}_t)$, and 
 $\langle M^i,M^j\rangle_t= \int_0^t
\bar{a}^{ij}(s,Y_s)ds$ for all  $i,j=1,\dots,d$. Moreover,  if $\Q^T$ is an extremal solution to the martingale problem $(MP)_Q$, then for any
martingale $(N_t)_{t\in [0,T]}$ with respect to  $\Q^T$ and $({\cal
G}_t)$ such that $N_0=0$, there exist predictable processes
$(h_t^j)_{t\in [0,T], j=1,\dots d}$ with  $\sum_{i,j=1} ^d
 \int_0^T h_s^i\bar{a}_{ij}(s,Y_s)h_s^j  ds<\infty, \  \Q^T  $ a.s.,   and such that  $(\int_0^t h_s
\cdot d M_s=\sum_{j=1}^d \int_0^t h_s^j d M_s^j )_{t\in [0,T]}$ is a modification of $(N_t)_{t\in
[0,T]}$. In particular, $(N_t)_{t\in [0,T]}$ has a continuous
modification.
\end{alem}

The main assertions in    part  b) of Lemma  \ref{DtTexplicit} are then  consequences of the next result.

\begin{alem}\label{Girsanovdensity}
Assume that $H1)$, $H2)_Q$, $H3)_Q$ and $H3)_P$ hold together. Suppose moreover
that $P_0 \ll Q_0$ and that $\Q^T$ is an extremal solution
to the martingale problem $(MP)_Q$. Recall that   $(t,x)\mapsto 
\nabla[\frac{p_{t}}{q_{t}}](x)$ is $ q_{t}(x)dx \ dt$  a.e. defined in $[0,T]\times \R^d \to \R^d$    by 
$
  \nabla \left[ \frac{p_{t}}{q_{t}}\right](x) : = \frac{p_{t}}{q_{t}}(x)   \nabla  \left[  \ln \frac{p_{t}}{q_{t}}\right](x) $.

\begin{itemize}
\item[i)]  With   $R$  the $({\cal G}_t)$-stopping time $R:=\inf\{s
\in [0,T]:D_s=0\}$, we have  $\Q^T-$a.s. that
\begin{equation*}
\begin{split} 
\forall t\in [0,T],\; &\int_0^t \nabla^*
\left[\frac{p_{T-s}}{q_{T-s}}\right](Y_s)  \bar{a}(s,Y_s)
\nabla \left[\frac{p_{T-s}}{q_{T-s}}\right](Y_s)\mathbf{1}_{s<R}  \ ds <\infty,  \mbox{ and }\\
\forall t\in [0,R),\; &\int_0^t \nabla^* \left[\ln
\frac{p_{T-s}}{q_{T-s}}\right](Y_s)  \bar{a}(s,Y_s) \nabla
\left[\ln \frac{p_{T-s}}{q_{T-s}}\right](Y_s) ds < \infty \mbox{ on }\{R>0\}.
\\
\end{split}
\end{equation*}
\item[ii)] The process  $(D_t)_{t\in
[0,T]}$ has a continuous version, denoted in the same way, such that 
\begin{equation*}
\begin{split}
\Q^T \, a.s,  \, \forall t\in [0,T],\;D_t= &  \frac{ p_T}{q_T}(Y_0) + \int_0^{t} \nabla
\left[\frac{p_{T-s}}{q_{T-s}}\right](Y_s) \mathbf{1}_{s<R} \cdot
dM_ s \\
= & \frac{ p_T}{ q_T}(Y_0) + \int_0^{t} \nabla
\left[\frac{p_{T-s}}{q_{T-s}}\right](Y_s) \mathbf{1}_{\{
\frac{p_{T-s}}{q_{T-s}}(Y_s)>0\}} \cdot dM_ s\ \\
\mbox{ and }\langle D \rangle_t=&\int_0^t \nabla^*
\left[\frac{p_{T-s}}{q_{T-s}}\right](Y_s)  \bar{a}(s,Y_s)
\nabla \left[\frac{p_{T-s}}{q_{T-s}}\right](Y_s)\mathbf{1}_{s<R} \ ds.
\end{split}
\end{equation*}
\end{itemize}
 
\end{alem}

\begin{adem}
  By Lemma \ref{reprmart}, the $\Q^T$-martingale $(D_t)_{t\in[0,T]}$ admits  the continuous version $D_0 +\sum_{j=1}^d\int_0^t h_s ^j
 d M_s^j$ still denoted by $D_t$ for simplicity. The martingale  representation property  and standard properties of stochastic integrals  moreover  imply that $D_t$  is determined by  the processes $\langle D,M^i\rangle=  \int_0 ^{\cdot} \bar{a}_{i j}(t,Y_t ) h^j_t  dt, \ i=1,\dots,d$. Consequently, $h_t$ can be replaced (leaving $D_t$ unchanged) by any predictable process $k_t$  such that for each $i$, $\int_0^{\cdot}\sum_{j=1} ^d h^j_t  \bar{a}_{i j}(t,Y_t ) dt=\int_0 ^{\cdot}\bar{a}_{i j}(t,Y_t ) k^j_t  dt$  $\Q^T$ a.s. (the fact  that $\int_0 ^{T}k^i_s \bar{a}_{i j}(s,Y_s ) k^j_s ds =\int_0 ^{T}\sum_{i,j=1} ^d h^j_s \bar{a}_{i j}(s,Y_s ) h^i_s ds<\infty $ $\Q^T$ a.s. then follows immediately).    Furthermore, since $D_t=D_{t\wedge R}$ by standard properties of nonnegative continuous martingales, we  may and shall assume that $\Q^T$ a.s.  $h_t=h_t \mathbf{1}_{t<R}=h_t\mathbf{1}_{D_t>0} $ for all $t\in [0,T]$.  Let us also notice that, by Fubini's Theorem, it  $\Q^T-$a.s. holds that $D_s=\frac{p_{T-s}}{q_{T-s}}(Y_s) $   (and 
 then   $\mathbf{1}_{\{ R> s\}} =\mathbf{1}_{\{\frac{p_{T-s}}{q_{T-s}}(Y_s)>0\}}$) for  a.e. $s\in [0,T]$. 
 
 Now, by our assumptions and Theorem \ref{TRMNS} a),  $\P^{T} \ll \Q^T$  are probability measures respectively solving 
 the martingale problems $(MP)_P$ and  $(MP)_Q$.  The processes
 $\int_0^{\cdot} \bar{b}_i^P(t,Y_t)dt $ and 
$ \int_0^{\cdot} \bar{b}_i^Q(t,Y_t)dt +  \int_0^{\cdot} (D_t)^{-1}h^j_t d \langle M^i,M^j\rangle_t$
 then are $\P^{T}-$ indistinguishable (see e.g.  Proposition 12.18 v) in \cite{Jac}). Using  these facts, the expression for $\langle M^i,M^j\rangle$ in Lemma  \ref{reprmart} and part ii) of Lemma \ref{gradlog}  we deduce first  that,  $\P^{T}-$a.s.,  
\begin{equation}\label{diffdrifts}
 \bar{b}_i^P(t,Y_t)- \bar{b}_i^Q(t,Y_t)= \bar{a}_{ij}(t,Y_t ) \left(  h^j_t \frac{q_{T-t}}{p_{T-t}}(Y_t)\right) = \bar{a}_{i \bullet}(t,Y_t)( K^p(T-t,Y_t)- K^q(T-t,Y_t))
\end{equation}
for   a.e. $t\in [0,T]$ and each $i$.  By part ii) of Lemma \ref{gradlog} we  then also get 
\begin{equation*}
 \int_0 ^{\cdot}  \bar{a}_{i j}(t,Y_t )h^j_t dt=  \int_0 ^{\cdot}\
\bar{a}_{i 	\bullet}(t,Y_t)( K^p(T-t,Y_t)- K^q(T-t,Y_t))\frac{p_{T-t}(Y_t)}{q_{T-t}(Y_t)}dt, \ i=1,\dots, d, 
\end{equation*}
$\P^{T}-$a.s., and then $\Q^T-$a.s. because of our assumption  on $h$. From these  identities and our previous discussion we deduce that we can choose $h_t=  \nabla \frac{p_{T-t}}{q_{T-t}}(Y_t)\mathbf{1}_{\{
\frac{p_{T-t}}{q_{T-t}}(Y_t)>0\}}=  \nabla \frac{p_{T-t}}{q_{T-t}}(Y_t) \mathbf{1}_{\{ R> t\}} $. This proves part ii). The first  property  of the process $\nabla  \frac{p_{T-t}}{q_{T-t}}(Y_t)$ in i) is thus consequence of the general properties of  $h$ in the representation formula for $D_t$.   The second assertion in i)  easily follows from the first one, taking into account  the  definitions of   $\nabla  \frac{p_{T-t}}{q_{T-t}}(Y_t) $ and $\nabla \ln \frac{p_{T-t}}{q_{T-t}}(Y_t)$  and the  properties of $D_t$.

\end{adem}

% \begin{arem}\label{equivclasse}

% We notice from the proof of Lemma \ref{Girsanovdensity} that  the function $\nabla  \frac{p_t}{q_t}(x)$  therein can be replaced by  any representant of the   equivalence class of the function $\frac{p_t}{q_t}(x)\left(K^p(t,x)- K^q(t,x)\right)$ under the relation $f(t,x)-g(t,x)\in Ker(a(t,x)), \ q_t(x)dx\ dt \  a.e.$.
% If $p_t$ and $q_t$ are of class  $C^1$  and $q_t$ does not vanish on $\R^d$, the true gradient $\nabla \frac{p_t}{q_t}$ is  equal to $\frac{p_t}{q_{t}} \left(\mathbf{1}_{p_t>0}\frac{\nabla p_t}{p_t} - \frac{\nabla q_t}{q_t}\right)$ (as the gradient of
% non-negative function vanishes at zeros of that  function) and is such a representant. 

% \end{arem}
 
\subsection{Proof of Theorem \ref{entprod}}

Since by Lemma \ref{DtTexplicit}, $(D_t)_{t\in[0,T]}$ is a continuous non-negative ${\mathbb Q}^T$-martingale and $U'_-$ is locally bounded on $(0,+\infty)$, $t\mapsto \int_0^t \left[U'_-(D_s)\right]^2 d\langle D\rangle_s$ is finite and continuous on $[0,T]$ when $R>T$   and finite and continuous on $[0,R)$ otherwise. In the latter case,
$\int_0^R \left[U'_-(D_s)\right]^2 d\langle D\rangle_s$ makes sense but is possibly infinite.
Define for any positive
integer $n$ the stopping time $$R_n:=\inf\left\{t\in [0,T\wedge R]:
D_t\leq \frac{1}{n}\mbox{ or }\int_0^t \left[U'_-(D_s)\right]^2 d\langle D\rangle_s \geq n\right\}.$$
For all $t\in [0,T]$, $\int_0^{t\wedge R_n}\left[U'_-(D_s)\right]^2 d\langle D\rangle_s\leq n$ and $\E\left(\int_0^{t\wedge R_n} U'_-(D_s)dD_s\right)=0$. Moreover $R_n\nearrow R$ as $n\to\infty$.

Let $t\in[0,T]$. By Tanaka's formula, \begin{equation}
\begin{split}
U(D_{t\wedge R_n }) =
 & U(D_0)+  \int_0^{t\wedge R_n} U'_-(D_s)dD_s+\frac{1}{2}\int_{(0,+\infty)}  L^r_{t\wedge R_n}(D)U''(dr).
\end{split}
\label{itoUstopm}\end{equation}
The
finiteness of $H_U(P_0\vert Q_0)$ implies that $(U(D_s))_{s\in[0,T]}$ is a uniformly integrable ${\mathbb Q}^T$-submartingale. Since the ${\mathbb Q}^T$-expectation of the stochastic integral is zero, one deduces
\begin{equation*}
\begin{split}
\E^T\left(U(D_{t\wedge R_n })\right) =  \E^T(U(D_0))+  \frac{1}{2} \E^T \left(\int_{(0,+\infty)} L^r_{t\wedge R_n}(D)U''(dr) \right).
\end{split}
\end{equation*}
When $n\to \infty$, since $U$ is continuous on $(0,+\infty)$ by convexity, $U(D_{t\wedge R_n })$ converges to $U(D_{t\wedge R})+\Delta U(0)1_{\{0<R\leq t\}}=U(D_t)+\Delta U(0)1_{\{0<R\leq t\}}$. Then, by uniform integrability, $\E^T(U(D_{t\wedge R_n }))$ converges to $\E^T(U(D_t))+\Delta U(0){\mathbb Q}^T(0<R\leq t)$. Dealing with the expectation of the integral in the right-hand-side above  by monotone convergence, we obtain
\begin{equation*}
\begin{split}
\E^T(U(D_t))=
\E^T(U(D_0))-\Delta U(0){\mathbb Q}^T(0<R\leq t)+\frac{1}{2}\E^T\left(\int_{(0,+\infty)} L^r_{t\wedge R}(D)U''(dr) \right).
\end{split}
\end{equation*}
Since according to Lemma \ref{DtTexplicit} b), $D$ is equal to zero on $[R,T]$, one can replace $t\wedge R$ by $t$ in the last expectation.
Replacing $t$ by $T-t$ in this equation, one gets \eqref{entderivloc}.
Moreover ${\mathbb Q}^T$ a.s., $\int_{(0,+\infty)}L^r_{t}(D)U''(dr) $ is the finite limit of the integral with respect to $U''(dr)$ in the right-hand-side of \eqref{itoUstopm} as $n\to\infty$. Since the left-hand side converges to $U(D_t)+\Delta U(0)1_{\{0<R\leq t\}}$ we deduce that the stochastic integral in the right-hand-side also has a finite limit. Hence $\int_0^{t\wedge R} \left[U'(D_s)\right]^2 d\langle D\rangle_s<+\infty$, $\int_0^{t\wedge R} U'(D_s)dD_s$ makes sense and \eqref{itoU} holds.
When $U$ is continuous on $[0,+\infty)$ and $C^2$ on $(0,+\infty)$, \eqref{itoUC2} follows by the  occupation times formula. In that case,  Lemma \ref{DtTexplicit} b) and \eqref{entderivloc} written for $t=0$ combined with the same arguments imply that
\begin{align*}H_U(P_0|Q_0)=& H_U(P_T|Q_T)\\&+\frac{1}{2}\E^T\left(\int_0^TU''(D_s)1_{\{s<R\}} \nabla^*\left[\frac{p_{T-s}}{q_{T-s}}\right](Y_s) \bar{a}(s,Y_s) \nabla
\left[\frac{p_{T-s}}{q_{T-s}}\right](Y_s)
ds\right).\end{align*} Since $Y_s$ admits the
density $q_{T-s}$ and for almost all $s\in [0,T)$,  $D_s=\frac{p_{T-s}}{q_{T-s}}(Y_s)$  and 
$\{R>s\}=\{\frac{p_{T-s}}{q_{T-s}}(Y_s)>0\}$,  \eqref{entderiv} follows by
the change of variables $s\mapsto T-s$.

\subsection{Proof of Corollary  \ref{TV}}\label{apptv} We  notice first that 
 \begin{equation}\label{finito}
  \forall \delta\in (0,1),\; \E^T \int_0^{T} \mathbf{1}_{|D_s -1|<\delta } \nabla^*
\left[\frac{p_{T-s}}{q_{T-s}}\right](Y_s) \bar{a}(s,Y_s)
\nabla \left[\frac{p_{T-s}}{q_{T-s}}\right](Y_s)ds <\infty.
\end{equation}
 Indeed, for $\delta\in (0,1)$, we can easily construct  a  $C^2$  convex function $\hat{U}$ on $\R$ such that $\forall r\in \R,\;0\leq \hat{U}(r)\leq |r-1|$    and     $\forall  r\in [1-\delta,1+\delta],\;\hat{U}''(r)\geq \alpha$ for some $\alpha >0$, so that the  integral  in \eqref{finito} is bounded thanks to \eqref{entderiv}  by $\frac{1}{\alpha} H_{\hat{U}}(P_0|Q_0) \leq  \frac{1}{\alpha} \|P_0-Q_0\|_{\rm TV}  $. 
  For $r\in\R$, since
  $$L_t^r(D)=2\left((D_t-r)^+-(D_0-r)^+ -\int_0^t \mathbf{1}_{D_s>r} dD_s\right),$$ 
by Doob's inequality we obtain
  $| \E^T ( L_t^r(D)-L_t^1(D))| \leq 4|r-1| + 2 \left( \E^T \int_0^t \mathbf{1}_{\{1\wedge r  < D_s\leq  r\vee 1\}} d\langle D\rangle _s\right)^{1/2}$. Hence, Lemma \ref{DtTexplicit} b) and  \eqref{finito} imply  that  $r\mapsto \E^T ( L_t^r(D))$ is continuous  (and finite) at  $r= 1$.   With the occupation times formula, one deduces that 
 \begin{equation*}
 \begin{split} 
 2 \E^T (L_t^1(D))   = & \lim_{\varepsilon\to 0}\frac{1}{\varepsilon}  \int_{1-\varepsilon}^{1+\varepsilon}
 \E^T (L^r_t (D))dr \\
 = & \lim_{\varepsilon\to 0}  \E^T\frac{1}{\varepsilon}  \int_0^{t} \mathbf{1}_{\{| D_s -1|<  \varepsilon \}} \nabla^* 
\left[\frac{p_{T-s}}{q_{T-s}}\right](Y_s)  \bar{a}(s,Y_s)
\nabla \left[\frac{p_{T-s}}{q_{T-s}}\right](Y_s)ds \\
=
& \lim_{\varepsilon\to 0} \int_0^{t}  \frac{1}{\varepsilon}  \int_{\{| \frac{p_{T-s}}{q_{T-s}} (x) -1|< \varepsilon\} }  \nabla^*
\left[\frac{p_{T-s}}{q_{T-s}}\right](x) \bar{a}(s,x)
\nabla \left[\frac{p_{T-s}}{q_{T-s}}\right](x) q_{T-s}(x) dx ds. \\
   \end{split}
   \end{equation*}
     Define now  the function $\varphi_{\varepsilon}(r):=\mathbf{1}_{[-\varepsilon,\varepsilon]}(r)r\varepsilon^{-1}  +\mathbf{1}_{(\varepsilon,\infty)}(r)-\mathbf{1}_{(-\infty,-\varepsilon)}(r) $.
  Since the function $\varepsilon \mapsto  \int_0^{t}    \int_{\{| \frac{p_{T-s}}{q_{T-s}} (x) -1|\leq \varepsilon \}}   q_{T-s}(x) dx ds $ is increasing and right continuous, we can chose  $\varepsilon_k \searrow 0$ a  sequence with $\int_0^{t}    \int_{\{| \frac{p_{T-s}}{q_{T-s}} (x) -1|=\varepsilon_k \}}   q_{T-s}(x) dx ds =0$ so that 
   \begin{align*} 
2\E^T (L_t^1(D))  & = 
   \lim_{k \to \infty} \int_0^{t}   \int_{\R^d} \nabla^*\left[ \varphi_{\varepsilon_k}\left( \frac{p_{T-s}}{q_{T-s}} -1\right)\right] (x) \bar{a}(s,x)
\nabla \left[\frac{p_{T-s}}{q_{T-s}}\right](x) q_{T-s}(x) dx ds \\
 &= -    \lim_{k \to \infty} \int_0^{t}   \int_{\R^d}  \varphi_{\varepsilon_k}\left( \frac{p_{T-s}}{q_{T-s}} -1\right) (x)  \nabla\cdot \left[ \bar{a}(s,x)
\nabla \left[\frac{p_{T-s}}{q_{T-s}}\right](x) q_{T-s}(x)\right] dx ds\\ 
&=-\int_0^{t}   \int_{\R^d}  \widetilde{sign} \left( \frac{p_{T-s}}{q_{T-s}} -1\right) (x)  \nabla\cdot \left[ \bar{a}(s,x)
\nabla \left[\frac{p_{T-s}}{q_{T-s}}\right](x) q_{T-s}(x)\right] dx ds 
   \end{align*}
where the last equality follows from the integrability assumption made on $\nabla\cdot \left[ \bar{a}(s,x)
\nabla \left[\frac{p_{T-s}}{q_{T-s}}\right](x) q_{T-s}(x)\right]$.  To justify the integration by parts at the second equality, we introduce functions $\phi_n\in C_0^{\infty}(\R^d)$  such that $\mathbf{1}_{B(0,r_n)} \leq \phi_n\leq\mathbf{1}_{B(0,2r_n)}$ and   $0\leq |\nabla \phi_n |\leq 2/r_n$, and functions $\varphi_{\varepsilon_k,m}:\R\to \R$ of class $C^1$ such that $\varphi_{\varepsilon_k,m}\to \varphi_{\varepsilon_k}$,  $|\varphi_{\varepsilon_k,m}| \leq |\varphi_{\varepsilon_k}|$ on $\R$  and $\varphi_{\varepsilon_k,m}'\to \varphi_{\varepsilon_k}' $, $|\varphi_{\varepsilon_k,m}' | \leq |\varphi_{\varepsilon_k}'|$ on $\R\backslash \{-\varepsilon_k,+\varepsilon_k\}$ as $m\to \infty$. Using the assumptions, \eqref{finito} and the choice of $\varepsilon_k$, we take the limits $n\to\infty$ then $m\to\infty$ by  dominated convergence in the equality
  \begin{align*}
    \int_{\R^d} \varphi_{\varepsilon_k,m}' &\left( \frac{p_{T-s}}{q_{T-s}} -1\right)(x) \nabla^*  \left[\frac{p_{T-s}}{q_{T-s}}\right](x)  a(T-s,x)
\nabla \left[\frac{p_{T-s}}{q_{T-s}}\right](x) q_{T-s}(x)  \phi_n(x)dx \\=&-\int_{\R^d}  \varphi_{\varepsilon_k,m}\left( \frac{p_{T-s}}{q_{T-s}} -1\right)(x) \nabla\cdot \left(a(T-s,x)
\nabla \left[\frac{p_{T-s}}{q_{T-s}}\right](x) q_{T-s}(x) \right) \phi_n(x)dx \\&- \int_{\R^d}  \varphi_{\varepsilon_k,m}\left( \frac{p_{T-s}}{q_{T-s}}  -1\right)(x)  \nabla^* \phi_n(x)   a(T-s,x)
\nabla \left[\frac{p_{T-s}}{q_{T-s}}\right](x) q_{T-s}(x)dx.
\end{align*}

\subsection{ Proof of Proposition \ref{vanish}}  \label{appvanish}
To check the Feller property, we introduce a continuous function $f:\R^d\to \R$ going to $0$ at infinity. Using It\^o's calculus and Gronwall's Lemma we check under the assumptions on the coefficients  that  the solution $X_t^x$  of   \eqref{diffhom} starting from $x\in \R^d$ satisfies $\E\left( (1+|X_t^x|^2)^{-1}\right) \leq C (1+|x|^2)^{-1}$ for some $C>0$.  Then, the inequality  
$$| \E(f(X_t^x))| \leq \sup_{|y|\leq A }  |f(y)| C \frac{ (1+A^2) }{(1+|x|^2)} +  \sup_{|y|> A }  |f(y)| $$ for all $A>0$ (following from the previous estimate and Markov's inequality) implies that $\E(f(X_t^x))\to 0$ when $x\to \infty$. Last, the continuity of $x\to\E(f(X^x_t))$ follows from the bound $\E(|X^x_t-X^y_t|^2)\leq C|x-y|^2$ and the uniform continuity and boundedness of $f$.

 By Theorem 1.3.8 \cite{Kun'}, since $(X_t)_{t\geq 0}$ is Feller  the tail sigma field is trivial as soon as  $\| P_t -Q_t \|_{TV}\to 0 \mbox{  as }t\to \infty$  for all pair of initial laws $P_0$ and $Q_0$. 
Since $\| P_t -Q_t \|_{TV}\leq \| P_t -p_\infty dx \|_{TV} + \| p_\infty dx- Q_t \|_{TV}$ and, by Theorem 2.1.3 p.162 \cite{boulhirsch}, the local uniform ellipticity assumption ensures that   $P_t$ admits a density with respect to the Lebesgue measure for all $t>0$, it is enough to show that $\| P_t -p_\infty dx \|_{TV}\to 0 \mbox{  as }t\to \infty$ when $P_0$ admits a density $p_0$ with respect to the Lebesgue measure.
 
 For $k\in\N^*$ 
consider the probability density $$p^k_0(x)=(p_0(x)\wedge kp_\infty(x))+p_\infty(x)\int_{p_0>k  p_{\infty} }(p_0(y)-kp_\infty(y))dy.$$

 Since $p_\infty$ is positive, on one hand we have $\lim_{k\to\infty}\|p_0-p_0^k\|_1=0$ and $p^k_0\leq (k+1)p_\infty$. On the other hand, the total variation distance between the marginal laws at time $t$ of the solutions to \eqref{diffhom} started from the initial densities $p_0$ and $p_0^k$ is not larger than $\|p_0-p_0^k\|_1$. Therefore we can moreover restrict ourselves to the case when $\frac{p_0}{p_\infty}$ is bounded. Then, $$\int_{\R^d}\left(\frac{p_0}{p_\infty}(x)-1\right)^2p_\infty(x)dx\leq 
\left(\int_{\R^d}\left(\frac{p_0}{p_\infty}(x)-1\right)^4p_\infty(x)dx\right)^{1/2}<+\infty.$$
We set $Q_0=p_\infty dx$. By Remarks \ref{condthmprincip} a) and \ref{condthmprincip} c), conditions $H1)$ , $H2)_{Q}$,  $H3)_{Q}$  and $H3)_{P}$ hold and for each $T>0$, $\Q^T$ is an extremal solution of the martingale problem $(MP)$. Applying Theorem \ref{entprod} respectively with $U(r)=(r-1)^4$ and $U(r)=(r-1)^2$, we get  that $t\mapsto \int_{\R^d}\left(\frac{p_t}{p_\infty}(x)-1\right)^2p_\infty(x)dx$ is non-increasing and  that \begin{equation}
   \sup_{t\geq 0}\int_{\R^d}\left(\frac{p_t}{p_\infty}(x)-1\right)^4p_\infty(x)dx+ \int_0^\infty\int_{\{\frac{p_{t}}{p_\infty}(x)>0\}}\left(\nabla^*\left[
\frac{p_{t}}{p_{\infty}}\right] a \nabla \left[\frac{p_{t}}{p_{\infty}}\right]\right)(x) p_{\infty}(x)dxdt<+\infty.\label{fini}
\end{equation}
Since $a$ is locally uniformly elliptic, the proof of Lemma \ref{gradlog} ensures that $dt$ a.e., the gradient $\nabla p_t$ (resp. $\nabla p_\infty$) of $p_t$ (resp. $p_\infty$) in the sense of distributions is a locally integrable function on $\R^d$ that vanishes a.e. on $\{x:p_t(x)=0\}$. Moreover, we can choose therein $K^p(t,x)=\mathbf{1}_{\{p_t(x)>0\}}\frac{\nabla p_t}{p_t}(x)$ and $K^q(t,x)=\frac{\nabla p_\infty}{p_\infty}(x)$. Then, in \eqref{fini}, $\nabla \left[
\frac{p_{t}}{p_{\infty}}\right]=\frac{\nabla p_t}{p_\infty}-\frac{p_t\nabla p_\infty}{p_\infty^2}$ is a.e. equal to $0$ when $\frac{p_t}{p_\infty}$ is equal to $0$ so that  the restriction of the spatial integration to $\{\frac{p_{t}}{p_\infty}(x)>0\}$ can be removed. Since $p_\infty$ is assumed to be locally Lipschitz continuous and bounded away from $0$, the function $\frac{1}{p_\infty}$ is locally bounded with a locally bounded distributional gradient equal to $-\frac{\nabla p_\infty}{p_\infty^2}$. We deduce that the gradient $\nabla\frac{p_t}{p_\infty}$ of $\frac{p_t}{p_\infty}$ in the sense of distributions is equal to $\frac{\nabla p_t}{p_\infty}(x)-\frac{p_t\nabla p_\infty}{p_\infty^2}$ and therefore to $\nabla \left[
\frac{p_{t}}{p_{\infty}}\right]$.

{F}rom the finiteness of the time-integral in \eqref{fini}, we deduce the existence of a sequence $(t_n)_n$ tending to $+\infty$ such that 
$\lim_{n\to \infty}\int_{\R^d}\left(\nabla^*
\frac{p_{t_n}}{p_{\infty}}a \nabla 
\frac{p_{t_n}}{p_{\infty}}\right)(x) p_{\infty}(x)dx=0$.
For $A>0$, writing the integral on $\R^d$ as the sum of the integrals on the ball $B(0,A)$ and its complementary $B(0,A)^c$, one has
\begin{align*}
   \int_{\R^d}&\left(\frac{p_{t_n}}{p_\infty}(x)-1\right)^2p_\infty(x)dx\\&\leq \int_{B(0,A)}\left(\frac{p_{t_n}}{p_\infty}(x)-\frac{\int_{B(0,A)} p_{t_n}(y)dy}{\int_{B(0,A)} p_{\infty}(y)dy}\right)^2p_\infty(x)dx+\frac{\left(\int_{B(0,A)}(p_{t_n}-p_\infty)(y)dy\right)^2}{\int_{B(0,A)}p_\infty(y)dy}\\&+
\left(\int_{B(0,A)^c}\left(\frac{p_{t_n}}{p_\infty}(x)-1\right)^4p_\infty(x)dx\int_{B(0,A)^c}p_\infty(x)dx\right)^{1/2}\\
&\leq \int_{B(0,A)}\left(\frac{p_{t_n}}{p_\infty}(x)-\frac{\int_{B(0,A)} \frac{p_{t_n}}{p_\infty}(y)dy}{\int_{B(0,A)} dy}\right)^2p_\infty(x)dx+\frac{\left(\int_{B(0,A)^c}(\frac{p_{t_n}}{p_\infty}(y)-1)p_\infty(y)dy\right)^2}{\int_{B(0,A)}p_\infty(y)dy}\\&+
\left(\int_{\R^d}\left(\frac{p_{0}}{p_\infty}(x)-1\right)^4p_\infty(x)dx\int_{B(0,A)^c}p_\infty(x)dx\right)^{1/2}.
\end{align*}Since $\left(\int_{B(0,A)^c}(\frac{p_{t_n}}{p_\infty}(y)-1)p_\infty(y)dy\right)^2 \leq \int_{\R^d}(\frac{p_{0}}{p_\infty}(y)-1)^2p_\infty(y)dy\int_{B(0,A)^c}p_\infty(y)dy$, the sum of the last two terms on the right-hand-side tends to $0$ uniformly in $n$ as $A\to\infty$. Using \eqref{localellip} and denoting by $C_A<+\infty$ the constant of the Poincaré-Wirtinger inequality satisfied by the Lebesgue measure on the ball $B(0,A)$, we check that the first term is smaller than $$C_A\frac{\sup_{B(0,A)}p_\infty}{\varepsilon_A\inf_{B(0,A)}p_\infty}\int_{\R^d}\left(\nabla^*
\frac{p_{t_n}}{p_{\infty}} a \nabla 
\frac{p_{t_n}}{p_{\infty}}\right)(x) p_{\infty}(x)dx,$$
which tends to $0$ as $n\to\infty$. Hence, $\lim_{n\to \infty}\int_{\R^d}\left(\frac{p_{t_n}}{p_\infty}(x)-1\right)^2p_\infty(x)dx=0$.
Since $ \| p_t -p_\infty \|_1^2\leq \int_{\R^d}\left(\frac{p_t}{p_\infty}(x)-1\right)^2p_\infty(x)dx$ where the right-hand-side is non-increasing with $t$, we conclude that $\lim_{t\to\infty}\|p_t-p_\infty\|_1=0$.

\subsection{Sufficient conditions for superquadratic potentials to satisfy $H1)''$}\label{proofsuperquad}

 \begin{alem}\label{superquad}

 Let     $b(x)=-\nabla V(x)$ for  a nonnegative $C^2$ potential $V$  in $\R^d$ satisfying \eqref{condv}, and $\sigma $ be any globally Lipschitz continuous choice of the square root of the identity $I_d$. 
Then, condition H1)'' holds  for the diffusion process $dX_t=\sigma(X_t)dW_t-\nabla V(X_t)dt$.
\end{alem}
\begin{adem}
Computing $d|X_t|^2$, we see that the first condition in \eqref{condv} prevents explosion for the SDE which has locally Lipschitz coefficients. Since for $c>0$, 
$$de^{cV(X_t)}=e^{cV(X_t)}\left(c\nabla^*V(X_t)\sigma(X_t)dW_t+\frac{c}{2}[\Delta V+(c-2)|\nabla V|^2](X_t)dt\right),$$
the second condition ensures that for $c$ small enough, $\E(e^{cV(X_t)})\leq e^{K(c)t}\E(e^{cV(X_0)})$ for some finite constant $K(c)$ only depending on $V$ and $c$. The third assumption ensures the existence of a finite constant $\tilde{K}(\frac{c}{T})$ only depending on $\frac{c}{T}$ and $V$ such that 
$$\E\left(\exp(4\int_0^T\sqrt{\partial_{ik}V\partial_{ik}V(X_t)}dt)\right)\leq \tilde{K}(\frac{c}{T})\E\left(\exp(\frac{c}{T}\int_0^TV(X_t)dt)\right).$$
By Jensen's inequality, we deduce that
$$\E\left(\exp(4\int_0^T\sqrt{\partial_{ik}V\partial_{ik}V(X_t)}dt)\right)\leq\frac{\tilde{K}(\frac{c}{T})}{T}\int_0^T\E(e^{cV(X_t)})dt\leq \tilde{K}(\frac{c}{T})e^{K(c)T}\E(e^{cV(X_0)}).$$
 \end{adem}

\section{Proofs of  the main results of  Section \ref{secdisfish} }\label{pdissipfish}

\subsection{Proof of Proposition \ref{dissipfish}}We will make use of  the
stochastic flow  defined by the two-parameter process $\xi_t(x)$  satisfying 
\begin{equation}\label{flow}
d\xi_t^i(x)= \sigma_{ik} (\xi_t(x))d \bar{W}^k_t
+ \bar{b}_i (\xi_t(x))dt, \quad (t,x) \in
[0,T)\times \R^d,  \  i=1, \dots d,
\end{equation}
and $ \xi_0(x)=x$,  noting that $\xi_t(Y_0)=Y_t$.  We shall also deal with the family  of continuous ${\cal G}_t-{\mathbb P}^{T}_{\infty}-$ local martingales $(D_t(x):t\in [0,T])_{x\in \R^d}$ defined by
\begin{equation}\label{Dtx}
dD_t(x) =\left[  \sigma_{ik}  \partial_i
\rho \right](t,\xi_t(x))
d\bar{W}_t^k \quad, 	\quad  D_0(x) =\frac{p_T}{p_{\infty}}(x)=\rho_0(x).
\end{equation}
According to Lemma \ref{DtTexplicit}, $D_t(Y_0)$ is equal to the process $D_t$ defined in \eqref{DtT}.
Writing $\nabla\rho_t(\xi_t(x))=(\nabla^*_x\xi_t(x))^{-1}\nabla_x[\rho_t(\xi_t(x))]$ we remark that, thanks to the It\^o product rule,  $d\nabla\rho_t(\xi_t(x))$ can be  obtained with  by computing   $d(\nabla_x\xi_t(x))^{-1}$ and $d\nabla_x[\rho_t(\xi_t(x))]$. Those computations are part of the contents  of the two next Lemmas:  

\begin{alem}\label{regflow}
The process $(t,x)\mapsto \xi_t(x)$ has a $\P^{T}_{\infty}$
a.s.  continuous version  such that the mapping $x\mapsto
\xi_t(x)$ is a global diffeomorphism of class $C^{1,\alpha}$ for
some $\alpha\in (0,1)$ and every $t\in [0,T]$. Moreover, we
have
\begin{equation}\label{derivflow}
d\partial_j\xi_t^i(x)= \partial_p \sigma_{ik} (t,\xi_t(x))\partial_j\xi_t^p(x) d
 \bar{W}^k_t
+\ \partial_p \bar{b}_i (t,\xi_t(x))\partial_j\xi_t^p(x)
dt,\quad (t,x) \in [0,T)\times \R^d
\end{equation}
 with $\partial_j \xi_0^i(x)=\delta_{ij}$. Finally, writing   $\nabla
\xi_t(x)=(\partial_j\xi_t^i(x))_{ij}$, it holds that
\begin{equation}\label{inversederivflow}
\begin{split}
d(\nabla \xi_t(x))^{-1}_{kl}= & -  (\nabla \xi_t
(x))^{-1}_{ki} [\partial_l \sigma_{ir}] (\xi_t(x)) d
\bar{W}^r_t -  \nabla \xi_t (x))^{-1}_{ki} [\partial_l
\bar{b}_{i}]
(\xi_t(x)) dt \\
& + (\nabla \xi_t (x))^{-1}_{ki} [\partial_m
\sigma_{ir}\partial_l \sigma_{mr}](\xi_t(x))dt,\ 
 \hspace{1cm} \quad (t,x) \in [0,T)\times \R^d.\\
\end{split}
\end{equation}
\end{alem}

\begin{adem} Under assumptions $H4)  $ and  $H5)_{p_\infty}$,  classic results of Kunita \cite{Kun'} (see Theorem 4.7.2) imply the
asserted regularity properties of the stochastic flow, as well as the
 $\P^{T}_{\infty}$ a.s. existence of the inverse matrix
$(\nabla \xi_t(x))^{-1}$  for all $(t,x)\in [0,T]\times
\R^d$. Since the smooth map $A\mapsto A^{-1}$, defined on non singular $d\times d$ matrices, has 
first and second derivatives respectively given by   the linear and bilinear operators $F\mapsto -A^{-1}FA^{-1}$ and
$(F,K)\mapsto A^{-1}F A^{-1} K A^{-1}+A^{-1}K A^{-1} F A^{-1}$ (where $F,K$ are generic square-matrices), we deduce  that for $A=(A_{ij})_{i,j=1\dots d}$,
$$\frac{\partial (A^{-1})_{kl}}{\partial
A_{ij}}=-A_{ki}^{-1}A_{jl}^{-1},\quad \mbox{ and } \quad
\frac{\partial^2 (A^{-1})_{kl}}{\partial A_{ij}\partial
A_{mn}}=A_{ki}^{-1}A_{jm}^{-1}A_{nl}^{-1}+A_{km}^{-1}A_{ni}^{-1}A_{jl}^{-1}$$
for all $k,l,i,j,m,n\in \{1,\dots,d\}$. Equation
\eqref{inversederivflow} follows  by  applying It\^o's formula to
each of the functions $A\mapsto (A^{-1})_{kl}$ and the
semimartingales $(\partial_j\xi_t^i(x))$, $i,j=1\dots d$.

\end{adem}

\begin{alem}\label{desintegD}
The process $D_t(x)$ has a  modification still denoted by $D_t(x)$
such that $\P^{T}_{\infty}$ a.s. the function $(t,x)\mapsto
D_t(x)$ is
 continuous and
 $x\mapsto D_t(x)$ is of class $C^1$ for each $t$.
 This modification is indistinguishable from   $(\rho_t(\xi_t(x)): (t,x)\in [0,T) \times \R^d)$
and   we have
\begin{equation}\label{Dxt'}
d\partial_k  D_t(x)=
 \partial_m \left[\sigma_{ir}
\partial_i \rho \right](t,\xi_t(x))
\partial_k \xi^m_t(x) d\bar{W}_t^r=d  \left[ \partial_m\rho(t,\xi_t(x))\partial_k \xi^m_t(x) \right]
\end{equation}
for all $(t,x)\in [0,T) \times \R^d$.
\end{alem}

\begin{adem}
Thanks to  the regularity of $x\mapsto
\xi_t(x)$  established in Lemma \ref{regflow} and assumptions $H5)_{p_\infty}$ and $H6)^T_{p_0}$, the statements  follow from  Theorem 3.3.3 of Kunita \cite{Kun'} (see also Exercise 3.1.5 therein).
\end{adem}

We can now  proceed to prove  Proposition \ref{dissipfish}. Evaluating expressions  \eqref{inversederivflow} and \eqref{Dxt'} in $x=Y_0$, we obtain using It\^o's product  rule that 
\begin{equation}\label{itoproduct}
\begin{split}
d \partial_l \rho_t(Y_t)= & 
\left[ \sigma_{kr}\partial_{lk}\rho \right](t,Y_t)
d\bar{W}^r_t - \left[\sigma_{kr}
\partial_{kj} \rho
\partial_l\sigma_{jr} +
\partial_k \rho \partial_l
 \bar{b}_k \right](t,Y_t)dt \\
 = & 
\left[ \sigma_{kr}\partial_{lk}\rho \right](t,Y_t)
d\bar{W}^r_t - \left[
 \frac{1}{2} \partial_{kj}  \rho \partial_l a_{kj} 
+
\partial_k \rho \partial_l
 \bar{b}_k \right](t,Y_t)dt. \\
 \end{split}
 \end{equation} 
{For} the remaining of the proof, the argument $(t,Y_t)$ will be omitted for notational simplicity.
By It\^o 's formula we get $
d\sigma_{li}= 
\left[\sigma_{mr}
\partial_{m} \sigma_{li}\right] 
d\bar{W}_t^r+\left[\bar{b}_m \partial_m
\sigma_{li}+
\frac{1}{2}a_{mk}\partial_{mk}\sigma_{li}
\right]dt.$ We then have 
\begin{equation*}
\begin{split}
d \left[ \sigma_{li}\partial_l \rho\right]= & \  \sigma_{li} d \partial_l \rho + \partial_l \rho d \sigma_{li} + d\langle
\partial_l \rho, \sigma_{li}\rangle \\
= &  \  \partial_{k }
\left[ \partial_l  \rho \sigma_{li} \right] \sigma_{kr}  d\bar{W}^r + \partial_l \rho \left[\bar{b}_m \partial_m
\sigma_{li}+
\frac{1}{2}a_{mk}\partial_{mk}\sigma_{li}
\right] -\sigma_{li} \left[\sigma_{kr}
\partial_{kj} \rho
\partial_l\sigma_{jr} +
\partial_k \rho \partial_l
 \bar{b}_k \right]  \\
& + \  a_{mk}  \partial_{lk}\rho  \partial_{m} \sigma_{li} \\
\end{split}
\end{equation*}
where we used  in the stochastic integral  the fact that $\partial_l \rho \sigma_{mr}\partial_m \sigma_{li} + \sigma_{li} \sigma_{kr} \partial_{lk} \rho=\partial_l \rho \sigma_{kr}\partial_k \sigma_{li} + \sigma_{li} \sigma_{kr} \partial_{lk} \rho=  \partial_{k }
\left[ \partial_l  \rho  \sigma_{li} \right] \sigma_{kr} $ . It follows that
\begin{equation*}
\begin{split}
d \left[\nabla^* \rho a \nabla \rho \right] = & \ d \left[ \sigma_{li}\partial_l \rho   \ 
 \sigma_{l'i}\partial_{l' }\rho \right] \\
=  &  \   2  \  \sigma_{l'i} \  \partial_{l'} \rho   \partial_{k }
\left[ \sigma_{li} \partial_l  \rho  \right] \sigma_{kr}   d\bar{W}^r+ 
 \ 2 \bigg\{   \left[\sigma_{l'i}\partial_{l'}\rho a_{mk}\partial_m\sigma_{li} \partial_{lk}\rho\right] \\
  & +  \sigma_{l' i} \partial_{l'}\rho \partial_l \rho \left[\bar{b}_m \partial_m
\sigma_{li}+
\frac{1}{2}a_{mk}\partial_{mk}\sigma_{li}
\right]   -  a_{ll'} \partial_{l'} \rho  \left[\sigma_{kr}
\partial_{kj} \rho
\partial_l\sigma_{jr} +
\partial_k \rho \partial_l
 \bar{b}_k \right] \bigg\} dt \\ 
 & + a_{kk'} \partial_k\left[\partial_l \rho \sigma_{li}\right] \partial_{k'}\left[\partial_{l' }\rho 
 \sigma_{l'i}\right] dt.
\end{split}
\end{equation*}    
 On the other hand, using \eqref{Dtx}   at $x=Y_0$  we have
$
d U''_{\delta}(\rho)=U^{(3)}_{\delta} (\rho)\sigma_{nr}\partial_n\rho \ d\bar{W}^r +\frac{1}{2} U^{(4)}_{\delta}(\rho)a_{nj}\partial_{n}\rho \partial_j \rho  \ dt$
 which combined with the previous expression yields
 \begin{equation}\label{disipinf}
\begin{split}
d \left[U''_{\delta} (\rho) \nabla^* \rho a \nabla \rho \right] = & 2 U''_{\delta} (\rho) 
\bigg\{ \left[\sigma_{l'i}\partial_{l'}\rho a_{mk}\partial_m\sigma_{li} \partial_{lk}\rho\right] +
  \sigma_{l' i} \partial_{l'}\rho \partial_l \rho \left[\bar{b}_m \partial_m
\sigma_{li}+
\frac{1}{2}a_{mk}\partial_{mk}\sigma_{li}
\right] \\ 
& \hspace{2cm}   -  a_{ll'} \partial_{l'} \rho  \left[\sigma_{kr}
\partial_{kj} \rho
\partial_l\sigma_{jr} +
\partial_k \rho \partial_l
 \bar{b}_k \right]  \bigg\}dt \  + \ d \hat{M}^{(\delta)} \\  
 & + \ U''_{\delta} (\rho)  a_{kk'}  \partial_k\left[\partial_l \rho \sigma_{li}\right] \partial_{k'}\left[\partial_{l' }\rho 
 \sigma_{l'i}\right] dt +  \frac{1}{2}U^{(4)}_{\delta} (\rho) \left|\nabla^* \rho a \nabla \rho \right|^2 dt \\
 & + \   2 U^{(3)}_{\delta} (\rho) \sigma_{l'i}\partial_{l'} \rho \partial_k \left[ \sigma_{li} \partial_l \rho \right] a_{jk}\partial_j \rho dt.
 \end{split}
 \end{equation}
 Equivalently, 
 \begin{equation*}
\begin{split}
 d \left[U''_{\delta}(\rho) \nabla^* \rho a \nabla \rho \right] = & 2 U''_{\delta}(\rho) 
\bigg\{ \partial_{l'}\rho \partial_l \rho \bigg[\frac{1}{4}(\partial_k\sigma_{lj}a_{km}\partial_{m}\sigma_{l'j}-\sigma_{ki}\partial_k\sigma_{lj}\sigma_{mj}\partial_{m}\sigma_{l'i})\\&\phantom{ 2 U''_{\delta}(\rho) 
\bigg\{ \partial_{l'}\rho \partial_l \rho \bigg[}+\frac{1}{2} \bar{b}_m \partial_m
a_{l l'}+\frac{1}{2}  \sigma_{l' i} a_{mk}  \partial_{mk}\sigma_{li}
-  a_{k l'}\partial_k
 \bar{b}_l \bigg] \\ 
& \hspace{2cm}   +\left[\sigma_{l'i}a_{mk}- \sigma_{ki}a_{ml'} \right] \partial_{l'} \rho
\partial_m\sigma_{li}\partial_{kl} \rho    \bigg\}dt \  +\ d \hat{M}^{(\delta)} +tr[\Lambda_\delta\Gamma]dt. 
 \end{split}
 \end{equation*}

\subsection{Proof of Theorem \ref{thdissipentro}}\label{pthdissipentro}
Let us check \eqref{minodissip2}.
Since $U'' $ is continuous and  non increasing in $(0,\infty)$ by Remark \ref{propertiesU}, one has $U''_{\delta}(r)\nearrow U''(r)$ for each $r>0$ as $\delta\to 0$. It is therefore  enough to obtain  (the integrated version of) inequality \eqref{minodissip2} with $U''_{\delta}$  instead of $U''$,  monotone convergence  allowing us to pass to the limit as $\delta\to 0$ on both sides.  
For  $0\leq r\leq  t < T$ we have by Proposition \ref{dissipfish}  that 
\begin{align}\label{finitevarpart2}
[U_{\delta} ''(\rho)& \nabla^* \rho a \nabla \rho](t ,Y_t)- [U''_{\delta} (\rho) \nabla^* \rho a \nabla \rho](r ,Y_{r}) \notag \\
&\geq  
\hat{M}_t^{(\delta)}-\hat{M}_r^{(\delta)}+2 \int_{r}^{t}   U_{\delta}''(\rho)\left[\sigma_{l'i}a_{mk}- \sigma_{ki}a_{ml'} \right] \partial_{l'} \rho
\partial_m\sigma_{li}\partial_{kl} \rho  
ds \ \notag \\
&+2\int_{r }^{t }  U''_{\delta} (\rho) \partial_{l'}\rho \partial_l \rho \bigg(\frac{1}{4}(\partial_k\sigma_{lj}a_{km}\partial_{m}\sigma_{l'j}-\sigma_{ki}\partial_k\sigma_{lj}\sigma_{mj}\partial_{m}\sigma_{l'i})\notag\\&\phantom{+2\int_{r }^{t }  U''_{\delta} (\rho) \partial_{l'}\rho \partial_l \rho \bigg(}+\frac{1}{2}\left[\bar{b}_m \partial_m
a_{ll'}  +\sigma_{l' i}
a_{mk}\partial_{mk}\sigma_{li}
\right]-  a_{m l'}\partial_m
 \bar{b}_l\bigg) ds.\end{align}
Since $ \partial_{ k l'}  \rho  U_{\delta}''(\rho)\left[\sigma_{l'i}a_{mk}- \sigma_{ki}a_{ml'} \right]=0 $ and  
$$ \partial_k (   U''_{\delta}(\rho) )  \partial_{l'} \rho  \left[\sigma_{l'i}a_{mk}- \sigma_{ki}a_{ml'} \right] = 
  U_{\delta}^{(3)}(\rho) \partial_k \rho   \partial_{l'}  \rho \left[\sigma_{l'i}a_{mk}- \sigma_{ki}a_{ml'} \right]  =0,$$one has
\begin{multline}\label{previntbypart}
U_{\delta}''(\rho)\left[\sigma_{l'i}a_{mk}- \sigma_{ki}a_{ml'} \right] \partial_{l'} \rho
\partial_m\sigma_{li}\partial_{kl} \rho  =    \frac{1}{p_{\infty}} \partial_k \left( \partial_{l}  \rho  \partial_{l'}  \rho  U_{\delta}''(\rho)\left[\sigma_{l'i}a_{mk}- \sigma_{ki}a_{ml'} \right]
\partial_m\sigma_{li}  \, p_{\infty}\right) \\
-  \frac{\partial_{l}  \rho  \partial_{l'}  \rho  U_{\delta}''(\rho)}{p_\infty}\partial_k \left( [a_{mk}\sigma_{l' i}- \sigma_{ki}a_{ml'}]\partial_m\sigma_{li}p_\infty \right).
  \end{multline}
  
Setting \begin{align*}
   \Sigma_{ll'}\stackrel{\rm def}{=}&\frac{1}{4}(\partial_k\sigma_{lj}a_{km}\partial_{m}\sigma_{l'j}-\sigma_{ki}\partial_k\sigma_{lj}\sigma_{mj}\partial_{m}\sigma_{l'i})+\frac{1}{2}\left[\bar{b}_m \partial_m
a_{ll'}  +\sigma_{l' i}
a_{mk}\partial_{mk}\sigma_{li}
\right]\\&-  a_{m l'}\partial_m
 \bar{b}_l-\frac{1}{p_\infty}\partial_k\left[\left(\frac{1}{2}a_{mk}\partial_ma_{ll'}- \sigma_{ki}a_{ml'}\partial_m\sigma_{li}\right)p_\infty \right]
\end{align*} we deduce that
\begin{align}\label{ineqdissipinf} 
&[U_{\delta} ''(\rho) \nabla^* \rho a \nabla \rho](t ,Y_t)- [U''_{\delta} (\rho) \nabla^* \rho a \nabla \rho](r ,Y_{r})\notag\\
&\geq\hat{M}_t^{(\delta)}-\hat{M}_r^{(\delta)}+2 \int_{r }^{t }  U''_{\delta}(\rho)    \Sigma_{ll'} \partial_{l'}\rho \partial_l \rho \, ds   +2\int_{r}^{t }   \frac{1}{p_{\infty}} \partial_k \left( \partial_{l}  \rho  \partial_{l'}  \rho  U_{\delta}''(\rho)\left[\sigma_{l'i}a_{mk}- \sigma_{ki}a_{ml'} \right]
\partial_m\sigma_{li}  \, p_{\infty}\right)  ds.
\end{align}

Using  \eqref{diffrev} and the identity  $\sigma_{ki} \partial_{k'}\sigma_{li} =  \partial_{k'}a_{kl} - \partial_{k'} \sigma_{ki} \sigma_{li}   $, one can check that
 \begin{equation}\label{thetalternative}
\begin{split}
\Theta_{ll'} = &  \frac{1}{2} \bar{b}_{k'} \partial_{k'} a_{l l' } + \frac{1}{2} ( a_{kl'}\partial_k\bar{b}_l+  a_{kl}\partial_k\bar{b}_{l'} )  + \frac{1}{4} a_{k'k}\partial_{k'k} a_{l l' }  -  \frac{1}{4} ( a_{k'k} \partial_{k'} \sigma_{li}\partial_k \sigma_{l'i}+\sigma_{ki}\partial_k\sigma_{lj}\sigma_{k'j}\partial_{k'}\sigma_{l'i})   \\
 & + \frac{1}{2}  \sigma_{ki} (\partial_{k'}\sigma_{li}a_{k'l'}+  \partial_{k'}\sigma_{l' i} a_{k'l})\partial_k\ln(p_\infty)  -\frac{1}{2}  a_{k'k}\partial_{k'} a_{l l' }\partial_k\ln(p_\infty)\\
 &+\frac{1}{2}\partial_k[   \sigma_{ki} (\partial_{k'}\sigma_{li}a_{k'l'}+  \partial_{k'}\sigma_{l' i} a_{k'l})  -a_{k'k}\partial_{k'} a_{l l' }  ]\\
 = & \frac{\Sigma_{ll'}+\Sigma_{l'l}}{2}
 \end{split}
 \end{equation} 
 and therefore, the second integral on the right-hand side of   \eqref{ineqdissipinf} rewrites as $2 \int_{r }^{t }  U''_{\delta}(\rho)    \Theta_{ll'} \partial_{l'}\rho \partial_l \rho \, ds $. 
 
Now, the quadratic variation of  $\hat{M}^{(\delta)}$  is bounded  above in  $ [0,T)$ by a constant times 
$$\int_0 ^{t}\left[ | U^{(3)}_{\delta}(\rho) |^2  | \nabla^*\rho  a\nabla\rho |^3 (Y_s)+ \left( U_{\delta}''(\rho) \right)^2   \nabla^*( \nabla^*\rho a\nabla\rho  ) a \nabla ( \nabla^*\rho a\nabla\rho) \right](Y_s) ds.$$ 
This fact and our assumptions imply that $\hat{M}^{\delta}$ is a martingale in $[0,T)$ for all $\delta>0$ sufficiently small.   Indeed, we have  from Remark \ref{propertiesU}  that  $U''_{\delta}(r)\leq U''(\delta)\wedge U''(r)$ and $|U^{(3)}_{\delta}(r)| \leq |U^{(3)} (\delta) |\wedge |U^{(3)}(r) |$ for all $r\geq 0$. Therefore (since $U''>0$) we have
 $U''_{\delta}(r)\leq (U''(r)\wedge 1) \1_{U''(\delta)\leq 1 }+ U''({\delta})(U''(r)/U''(\delta))\wedge 1) \1_{U''(\delta)> 1 }$ whence  $U''_{\delta}(r)\leq (U''(\delta)+ 1)(  U''(r)\wedge 1)$.  As $U^{(3)}$ is non decreasing and non positive, either $| U^{(3)}(\delta)|\not =0$ for all $\delta$ sufficiently  small, in which case we similarly get  $|U^{(3)}_{\delta}(r)|\leq (|U^{(3)}(\delta)|+ 1)(  |U^{(3)}(r)| \wedge 1)$, or otherwise $U^{(3)}_{\delta}$  identically vanishes for all $\delta$.  Assumption $H6')_{p_{\infty}}$  and the previous then ensure that   $\langle M^{(\delta)}\rangle_t $ has finite expectation for $t\in [0, T)$.
 
In order to conclude that inequality \eqref{minodissip2} holds for the function $U_{\delta}$ , noting that $\nabla \rho_t$ vanishes on $\{\rho_t=0\}$,  it is enough to show that the last integral  in \eqref{ineqdissipinf}  has (well defined) null expectation.  Using \eqref{previntbypart} and Assumption  $H6')_{p_{\infty}}$  we 
obtain (with the same estimation for  $U''_{\delta}(r)$ as before) that 
\begin{multline}\label{expectwelldefined}
{\mathbb E}^{T}_{\infty} \int_{r}^{t } \left|  \frac{1}{p_{\infty}}  \partial_k \left( \partial_{l}  \rho  \partial_{l'}  \rho  U_{\delta}''(\rho)\left[\sigma_{l'i}a_{mk}- \sigma_{ki}a_{ml'} \right]
\partial_m\sigma_{li}  \, p_{\infty}\right) \right|  (Y_s) ds \\
 =\int_{r}^{t } \int_{\R^d} \big| \partial_k \left( \partial_{l}  \rho  \partial_{l'}  \rho  U_{\delta}''(\rho)\left[\sigma_{l'i}a_{mk}- \sigma_{ki}a_{ml'} \right]
\partial_m\sigma_{li}  \, p_{\infty}\right) \big| dx  ds  <\infty\\
\end{multline}
which shows that the expectation of the last  term in \eqref{ineqdissipinf}  is well defined. Moreover, the (everywhere defined) spatial divergence of  $g(s,x):=\partial_{l}  \rho_s  \partial_{l'}  \rho_s  U_{\delta}''(\rho_s)\left[\sigma_{l'i}a_{m\bullet}- \sigma_{\bullet i}a_{ml'} \right]
\partial_m\sigma_{li}  \, p_{\infty}$ is $L^1(dx,\R^d)$ for a.e. $s$. For such $s$ and $\phi_n\in C_0^{\infty}(\R^d)$  satisfying $0\leq \phi_n\leq 1$,  $0\leq |\nabla \phi_n |\leq 1$, $\phi_n(x)=1$ for $x\in B(0,n)$ and  $\phi_n(x)=0$ for $x\in B(0,2n)^c$, we have
\begin{align*}
   0=\int_{\R^d}\nabla.(\phi_n(x)g(s,x))dx=\int_{\R^d}\phi_n(x)\nabla.g(s,x)dx+\int_{\R^d}\nabla\phi_n(x).g(s,x)dx.
\end{align*}
Since by Lebesgue's theorem, the second term of the right-hand-side tends to $0$ as $n\to\infty$, the limit $\int_{\R^d}\nabla.g(s,x)dx$ of the first term is equal to $0$.

\section{Dissipation of the Fisher information : comparison with the computations and results in \cite{Arnoldcarlenju}}\label{compacj}

 In this section we compare our computations and results with those in \cite{Arnoldcarlenju}.

 The form of the term $tr(\Lambda_{\delta}\Gamma)$ in Proposition \ref{dissipfish} is inspired from the term $tr(\mathbf{X}\mathbf{Y})$ in \cite{Arnoldcarlenju} pp 163-164 where $\mathbf{X}=2\Lambda_{\delta}$. One has 
  \begin{equation*}
  \begin{split}
\Gamma_{12}= \ (\nabla^*  \rho \  a)_j  \   \partial_j (\sigma_{ki}\partial_k \rho)\sigma_{li}\partial_l \rho=  &
 \frac{1}{2}(\nabla^*  \rho \  a)_j  \left[  \partial_j (\sigma_{ki}\partial_k \rho)\sigma_{li}\partial_l \rho +
  \partial_j (\sigma_{li}\partial_l \rho) \sigma_{ki}\partial_k \rho \right]   \\  =  & 
\frac{1}{2}(\nabla^*  \rho \  a)_j  \partial_j\left[  \partial_l \rho a_{kl}  \partial_k \rho\right] = 
\frac{1}{2}(\nabla^*  \rho \  a)\nabla ( \nabla^*  \rho a \nabla \rho )\\
 \end{split}
  \end{equation*}
  which,  with $\frac{\partial v}{\partial x}:=(\partial_j  v_i)_{i,j}$ denoting the Jacobian matrix of  vector field $v$,   equals
    \begin{equation*}
  \begin{split}
\frac{1}{2}(\nabla^*  \rho \  a)_j  \partial_j\left[ \partial_k \rho  a_{kl}  \partial_l \rho\right] =  &  \frac{1}{2}(\nabla^*  \rho \  a)_j \left(    \partial_{kj} \rho \  a_{kl} \   \partial_l \rho  +\partial_j\left[    a_{kl} \ \partial_l    \rho \right]   \partial_k \rho  \right)  \\ 
 =& \frac{1}{2}\nabla^*  \rho \  a  \frac{\partial (  \nabla   \rho)}{\partial x}  a \nabla \rho + \frac{1}{2}
  \nabla^*  \rho \  a   \frac{\partial ( a  \nabla   \rho)}{\partial x}^*  \nabla \rho \\
 \end{split}
  \end{equation*}
  and corresponds to $4\mathbf{Y}_{12}$ in \cite{Arnoldcarlenju} p. 164 (noting that in their notation, $\mathbf{D}(x)=a(x) / 2$). Similarly, $\Gamma_{22}=4\mathbf{Y}_{22}$.  However $\Gamma_{11}$ cannot  in general be identified with $4\mathbf{Y}_{11}$. For instance, in the case of scalar diffusion $\mathbf{D}(x)=a(x) / 2=D(x)I_{d}$ for some real valued function $D$, the 
 term $\Gamma_{11}(x)$ above when written in terms of ${D}$ reads 
 $$ \frac{1}{2}| \nabla D|^2 | \nabla \rho |^2 +\frac{1}{2}(\nabla D.\nabla\rho)^2+4 D\partial_{j} D\partial_ i \rho \partial_{ij} \rho + 4 D^2\sum_{ij}(\partial_{ij}\rho)^2$$ 
for the choice $\sigma(x)=\sqrt{2D(x)}I_d$,  whereas  
 $$4\mathbf{Y}_{11}=4\left(D^2 \sum_{ij}(\partial_{ij}\rho)^2 +\left(\frac{d}{4}-\frac{1}{2}\right) (\nabla \rho \cdot \nabla D)^2 +2D\partial_j D \partial_i \rho\partial_{ij}\rho   -D (\nabla \rho \cdot \nabla D)\triangle \rho +\frac{1}{2} |\nabla D|^2 |\nabla \rho|^2\right).$$
 Moreover,   our term  $\Gamma_{11}$ is   non-intrinsic, in  the sense that it cannot in general   be written in terms of the diffusion matrix  $a$ only (without making explicit use of $\sigma$), contrary to the term  $\mathbf{Y}_{11}$  in the matrix  of  \cite{Arnoldcarlenju}.

We will next  check that  the criterion in  \cite{Arnoldcarlenju} can also be derived from the computations in   Proposition \ref{dissipfish}   in case $a$ is non singular, which amounts to make an alternative choice in the expression  for $ d \left[U''_{\delta} (\rho) \nabla^* \rho a \nabla \rho \right]$ of the quantities in the roles of the coefficient   $\Gamma_{11}$ and  of the term  $\bar{\theta}$.   This will also allow us to compare and combine both criteria. 

Recall first that the matrix $\mathbf{D}(x) $ in  \cite{Arnoldcarlenju} equals half of  our matrix $a(x)$, and notice that our forward drift term writes in their notation 
$b=-\mathbf{D} \nabla \phi -\mathbf{D} F +\nabla. \mathbf{D}$,
where $(\nabla.  \mathbf{D} )_i=\partial_j  \mathbf{D}_{ij}$,  $e^{-\phi}=p_{\infty}$ is the invariant density,     and $F$ a  is vector field  satisfying $\nabla. ( \mathbf{D} F e^{-\phi})=0$. Thus, 
$\bar{b}=a \nabla \ln p_{\infty}   +\nabla. \, a -b = -\mathbf{D} \nabla \phi + \mathbf{D} F +\nabla. \mathbf{D}. $
%(in particular, the symmetric  and non symmetric parts of the drift are respectively given  by $\frac{b+\bar{b}}{2}=\frac{1}{2}\frac{\nabla. (a p_{\infty})}{p_{\infty}}=-  \mathbf{D}\nabla \phi +\nabla. \mathbf{D}$ and $ \frac{b-\bar{b}}{2}=-  \frac{1}{2}a F = -   \mathbf{D} F $). 
 
The factor of $ U''_{\delta}(\rho)$ in \eqref{disipinf} takes the intrinsic form
\begin{align*}
&a_{kk'}\left[\partial_{kl}\rho \sigma_{li} \partial_{k'l'}\rho  \sigma_{l'i} + \partial_{kl}\rho \sigma_{li}\partial_{l'}\rho \partial_{k'} \sigma_{l'i}+ \partial_{l} \rho \partial_{k} \sigma_{li}\partial_{k' l'}\rho \sigma_{l'i} + \partial_l \rho \partial_k \sigma_{li} \partial_{l'} \rho \partial_{k'} \sigma_{l'i}\right] \\  & + 2 \sigma_{l'i}\partial_{l'}\rho a_{k'k} \partial_{k'}\sigma_{li} \partial_{lk} \rho   + \partial_l\rho  \partial_{l'} \rho a_{kk'} \partial_{k'k} \sigma_{li} \sigma_{l'i}-2 a_{ll'} \partial_{l'}\rho\sigma_{kr}\partial_{kk'}\rho \partial_l \sigma_{k'r} +\bar{b}_m\partial_m a_{ll'}\partial_l\rho\partial_{l'}\rho-2a_{ll'}\partial_{l'}\rho\partial_k\rho\partial_l\bar{b}_k\\&=
 a_{kk'} \left[ \partial_{kl}\rho   \partial_{k'l'}\rho a_{ll'}+ 2 \partial_{kl}\rho  \partial_{l'}\rho  \partial_{k'} a_{ll'}\right] \\&+ \frac{1}{2} a_{kk'} \partial_l\rho  \partial_{l'} \rho \partial_{kk'} a_{ll'}- a_{ll'} \partial_{l'}\rho \partial_{kk'} \rho \partial_l a_{kk'}+\bar{b}_m\partial_m a_{ll'}\partial_l\rho\partial_{l'}\rho-2a_{ll'}\partial_{l'}\rho\partial_k\rho\partial_l\bar{b}_k, 
\end{align*}
where to  the second and third terms in the bracket on the left-hand side,  brought together, we have added the first term after the bracket, and moreover  the fourth term in the bracket  on the left-hand side was added to the the second term outside the bracket.  
Hence, writing     $$Q_1:= - a_{ll'} \partial_{l'}\rho \partial_{kk'} \rho \partial_l a_{kk'}+\bar{b}_m\partial_m a_{ll'}\partial_l\rho\partial_{l'}\rho-2a_{ll'}\partial_{l'}\rho\partial_k\rho\partial_l\bar{b}_k, $$ 
$$Q_2:=
 a_{kk'} \left[ \partial_{kl}\rho   \partial_{k'l'}\rho a_{ll'}+ 2 \partial_{kl}\rho  \partial_{l'}\rho  \partial_{k'} a_{ll'}\right] + \frac{1}{2} a_{kk'} \partial_l\rho  \partial_{l'} \rho \partial_{kk'} a_{ll'}, $$ 
   and using the last expression  for $\Gamma_{12}$  above, we can write  \begin{equation}
 \label{propdissiprewrites}
 \begin{split}
 \frac{1}{2}d& \left[  U''_{\delta}(\rho) \nabla^* \rho a \nabla \rho \right] =    \frac{1}{2} d \hat{M}^{(\delta)} + 
   \frac{U''_{\delta}(\rho)}{2}   (Q_1+Q_2)  \, dt 
    +  \frac{U^{(4)}_{\delta}(\rho)}{4} \left|\nabla^*  \rho a \nabla \rho \right|^2 dt  \\ 
    &\phantom{\left[  U''_{\delta}(\rho) \nabla^* \rho a \nabla \rho \right] =} +  \frac{U^{(3)}_{\delta}(\rho)}{2}  \left( 
    \nabla^*  \rho \   a   \frac{\partial (  \nabla   \rho)}{\partial x}   a  \nabla \rho + 
  \nabla^*  \rho \   a    \frac{\partial (  a   \nabla   \rho)}{\partial x}^*  \nabla \rho \right)dt.  \\ 
%  \begin{align*}
 = &  \frac{1}{2} d \hat{M}^{(\delta)}+  \left[ \frac{U^{(4)}_{\delta}(\rho) }{4}  \left|\nabla^*  \rho a  \nabla \rho \right|^2   +  \frac{U^{(3)}_{\delta}(\rho) }{4} \left( 
    \nabla^*  \rho \   a   \frac{\partial (  \nabla   \rho)}{\partial x}   a  \nabla \rho + 
  \nabla^*  \rho \   a    \frac{\partial (  a   \nabla   \rho)}{\partial x}^*  \nabla \rho \right)\right] dt    \\
   &  +\left[ \frac{U''_{\delta}(\rho)}{2}   Q_1 \,\right]  dt  +  \left[  \frac{U''_{\delta}(\rho) }{2}    Q_2 +   \frac{U^{(3)}_{\delta}(\rho) }{4} \left( 
    \nabla^*  \rho \   a   \frac{\partial (  \nabla   \rho)}{\partial x}   a  \nabla \rho + 
  \nabla^*  \rho \   a    \frac{\partial (  a   \nabla   \rho)}{\partial x}^*  \nabla \rho \right)\right] dt   .\\ 
 %   \end{align*}
  \end{split}
\end{equation}
The latter identity yields  the expression  for the dissipation of entropy dissipation computed in \cite{Arnoldcarlenju}. Indeed, 
denoting  respectively by  $J_1$, $J_2$  and $ J_3$ the expectations of the first, second and third terms in square brackets in the right-hand side,  we observe that  $J_1$  is, up to time reversal $t\mapsto T-t$,  exactly equal to the term $\tilde{R}_1$ on top of p. 162 in  \cite{Arnoldcarlenju}.  Starting from the last expression of $T_3$ p. 160 and the definition (2.23) of $\tilde{R}_2$ and $T_4$ %\cite{Arnoldcarlenju} 
and replacing $\mathbf{D} F$ by its expression $\bar{b}-\frac{1}{2}(a\nabla\ln(p_\infty)+\nabla.a)$  in our notation, we get  that $\tilde{R}_2+T_3+T_4$ is equal to 
\begin{align*}
   &\int_{\R^d}\frac{U''_{\delta}(\rho)Q_1}{2}p_\infty-\int_{\R^d}\left[\frac{U''_{\delta}(\rho)}{4}\left(\partial_i\rho \partial_j a_{ik}\partial_k\rho+2\partial_{ij}\rho a_{ik}\partial_k\rho\right)\right]\times\partial_l(a_{lj}p_\infty)-\int_{\R^d}\frac{U''_{\delta}(\rho)}{2}\left[a_{lj}\partial_{ijl}\rho\partial_k\rho a_{ki}\right]p_\infty
\end{align*}
up to time reversal. 
The first term corresponds to $J_2$. 
Integrating by parts the second term to get rid of the derivative with respect to the $l$-th coordinate in the second factor, one checks that its sum with the last one is equal to $J_3$. Hence, up to time reversal, we have  $J_1+J_2+J_3=(\tilde{R}_1+T_3)+ (\tilde{R}_2+T_4) $ which is the expression  for the dissipation of entropy dissipation computed in \cite{Arnoldcarlenju} p. 160.

In order to recover the  Bakry Emery criterion in \cite{Arnoldcarlenju}, we rewrite  $Q_1+Q_2=K_1(\rho)+K_2(\rho)$ where $$K_1(\rho):= \bar{b}_m\partial_m a_{ll'}\partial_l\rho\partial_{l'}\rho-2a_{ll'}\partial_{l'}\rho\partial_k\rho\partial_l\bar{b}_k+ \frac{1}{2} a_{kk'} \partial_l\rho  \partial_{l'} \rho \partial_{kk'} a_{ll'} $$ 
and
 $$K_2(\rho):=  a_{kk'}  \partial_{kl}\rho  \, a_{ll'}\partial_{k'l'}\rho     + 2a_{kk'}   \partial_{kl}\rho  \partial_{l'}\rho  \partial_{k'} a_{ll'} -a_{k' l'} \partial_{l'}\rho \partial_{kl} \rho \partial_{k'} a_{kl }.$$
 When $a$ is non singular, introducing  $G_{jk}(\rho)= \partial_{l'} \rho a_{k'l'} \partial_{k'}a_{jk}$ and $H_{lj}(\rho)=\partial_j a_{ll'} \partial_{l'} \rho$  we can  write
 \begin{equation*}
 \begin{split}
 K_2(\rho)  = & tr \left[(a\nabla^2\rho)^2 + 2 H(\rho)  \, a\nabla^2 \rho -G(\rho) \nabla^2\rho \right] \\
 = &  tr \left[(a\nabla^2\rho)^2 +  H(\rho) a\nabla^2 \rho   +  a H(\rho) ^* \nabla^2 \rho   -G(\rho) a^{-1} a \nabla^2\rho \right]  \\ 
= &  tr \bigg[(a\nabla^2\rho)^2   +  \frac{1}{2}( H(\rho) a\nabla^2 \rho   +  a H(\rho)^* \nabla^2 \rho   -G(\rho) a^{-1} a \nabla^2\rho ) \\
&  +\frac{1}{2}(a\nabla^2 \rho H(\rho) \, +  a\nabla^2 \rho a H(\rho)^*a^{-1} - a \nabla^2\rho G(\rho) a^{-1} )  \bigg]  \\
\end{split}
\end{equation*}
where we have  used the cyclicity of the  trace and its invariance by transposition. Following  \cite{Arnoldcarlenju}, we complete the trace of a squared sum of matrices  to get 
$$K_2(\rho) = tr \left[ a\nabla^2\rho + \frac{1}{2}( H(\rho) +   a H(\rho)^* a^{-1} -G(\rho) a^{-1})\right]^2 -\frac{1}{4}tr \left[ H(\rho) +   a H(\rho)^* a^{-1} -G(\rho) a^{-1}\right]^2.   $$
The finite variation part on  the right-hand side of the   first line in \eqref{propdissiprewrites} therefore rewrites
 \begin{equation}\label{finitevarpart}
  \begin{split}
  &  \frac{U''_{\delta}(\rho)}{2}( K_1(\rho) -\frac{1}{4}tr \left[ H(\rho) +   a H(\rho)^* a^{-1} -G(\rho) a^{-1}\right]^2)\\
&   
 +  \frac{U''_{\delta}(\rho)}{2}   tr \left[ a\nabla^2\rho + \frac{1}{2}( H(\rho) +   a H(\rho)^* a^{-1} -G(\rho) a^{-1})\right]^2  
     \, dt 
    \\ 
&     +   \frac{U^{(3)}_{\delta}(\rho)}{2}  \left( 
    \nabla^*  \rho \   a   \frac{\partial (  \nabla   \rho)}{\partial x}   a  \nabla \rho + 
  \nabla^*  \rho \   a    \frac{\partial (  a   \nabla   \rho)}{\partial x}^*  \nabla \rho \right)dt  \\ 
  & +  \frac{U^{(4)}_{\delta}(\rho)}{4} \left|\nabla^*  \rho a \nabla \rho \right|^2 dt . \\
\end{split}
\end{equation}
 The sum of the second, third and fourth lines correspond to the matrix product $\mathbf{X}\mathbf{Y}$ in  \cite{Arnoldcarlenju} and is shown to be nonnegative  in p. 164 therein.   We can then check that  for a smooth  function $v:\R^d\to \R$, the term $\frac{1}{2}( K_1(v) -\frac{1}{4}tr \left[ H(v) +   a H(v)^*  a^{-1} -G(v) a^{-1}\right]^2)$ is twice the expression on the left-hand side of   the inequality  (2.13) in p. 158 of \cite{Arnoldcarlenju} (with $\nabla v$ corresponding to their vector field ``$U$''). Consequently,  their   Bakry Emery criterion  (2.13)  corresponds, in our notation, to imposing the condition  
 
 $\exists \lambda>0$  such that  for all smooth function $v:\R^d \to \R$ and all $x\in \R^d$: 
 $$  \frac{1}{2}( K_1(v) -\frac{1}{4}tr \left[ H(v) +   a H(v)^* a^{-1} -G(v) a^{-1}\right]^2)(x) \geq  \lambda  \nabla v^* a \nabla v (x)  , 
 $$
 which implies exponential  convergence at rate $2\lambda$ of the $U-$Fisher information and the $U-$ relative entropy.
  
 We may combine this criterion with ours by introducing some $C^1$ function $\alpha: \R^d\to [0,1]$ and writing the finite variation part on the  right-hand side of the   first line in \eqref{propdissiprewrites}  as $(1-\alpha) $ multiplied by 
   the expression \eqref{finitevarpart}, plus $\frac{1}{2}\alpha $ multiplied by the finite variation part in the right-hand side of \eqref{finitevarpart2}. Because of the integration by parts performed in the proof of Theorem \ref{thdissipentro}, the mixed criterion involves the derivatives of $\alpha$. Let 
   $$\Theta^{\alpha}_{ll'}:=\alpha\Theta_{ll'}-\frac{1}{2}\partial_k \alpha \left( \left[\sigma_{l'i}a_{mk}- \sigma_{ki}a_{ml'} \right]
\partial_m\sigma_{li} +   \left[\sigma_{li}a_{mk}- \sigma_{ki}a_{ml} \right]
\partial_m\sigma_{l'i}  \right).$$
This ultimate mixed criterion writes

$\exists \lambda>0$  such that  for all smooth function $v:\R^d \to \R$ and all $x\in \R^d$: 
 $$\nabla v^* \Theta^{\alpha} \nabla v (x) +(1-\alpha(x)) \left(  \frac{1}{2}( K_1(v) -\frac{1}{4}tr \left[ H(v) +   a H(v)^* a^{-1} -G(v) a^{-1}\right]^2)(x) \right) \geq  \lambda  \nabla v^* a \nabla v (x)  
 $$
and also implies exponential  convergence at rate $2\lambda$ of the $U-$Fisher information and the $U-$ relative entropy.

\newpage

\end{document}